\documentclass[aap,MSNbibl,seceqn,dvips]{arximspdf}
\usepackage{url,breakurl}
\doi{10.1214/13-AAP992} 
\volume{25}
\issue{1}
\pubyear{2015}
\firstpage{177}
\lastpage{210}

\makeatletter
\newproclaim{rems}{Remarks}

\newcommand{\rrvert}{\vert}
\newcommand{\llvert}{\vert}
\def\cal{\mathcal}

\newcommand{\Y}{{\cal Y}}

\newcommand{\X}{{\cal X}}
\newcommand{\0}{\mathbf{0}}

\newtheorem{theo}{Theorem}[section]
\newtheorem{coro}{Corollary}[section]

\newtheorem{lemm}{Lemma}[section]
\newproclaim{defn}{Definition}[section]
\newproclaim{assu}{Assumption}[section]
\newproclaim{rema}{Remark}[section]
\def\HH{\mathbb{H}}
\def\E{\mathbb{E}}
\def\0{\mathbf{0}}
\def\Z{\mathbb{Z}}
\def\R{\mathbb{R}}
\def\B{{ B}}
\def\H{{\cal H}}
\renewcommand{\E}{\mathbb E}
\newcommand{\C}{{\cal C}}
\newcommand{\MM}{\mathbb{M}}

\newcommand{\tT}{\tilde{T}}
\newcommand{\tQ}{\tilde{Q}}

\newcommand{\M}{{\cal M}}
\newcommand{\Q}{{\cal Q}}
\renewcommand{\P}{{{\cal P}}}

\newcommand{\A}{{\cal A}}
\newcommand{\Cov}{\operatorname{Cov}}
\newcommand{\Var}{\operatorname{Var}}

\newcommand{\Vol}{\mathrm{Vol}}
\newcommand{\diam}{\operatorname{diam}}
\newcommand{\F}{{\cal F}}
\newcommand{\V}{{\cal V}}

\def\R{\mathbb{R}}

\def\al{{\alpha}}
\def\A{{\cal A}}
\def\la{{\lambda}}
\def\ka{{\kappa}}

\newcommand{\eqref}[1]{(\ref{#1})}
\makeatother

\begin{document}
\begin{frontmatter}

\title{Surface order scaling in stochastic geometry\thanksref{T1}}
\runtitle{Surface order scaling in stochastic geometry}
\thankstext{T1}{Supported in part by NSF Grant DMS-11-06619.}

\begin{aug}
\author[A]{\fnms{J. E.} \snm{Yukich}\corref{}\ead[label=e1]{joseph.yukich@lehigh.edu}}
\runauthor{J. E. Yukich}
\affiliation{Lehigh University}
\address[A]{Department of Mathematics\\
Lehigh University\\
Bethlehem, Pennsylvania 18015\\
USA\\
\printead{e1}} 
\end{aug}

\received{\smonth{4} \syear{2013}}
\revised{\smonth{11} \syear{2013}}

%
\begin{abstract} 
Let $\mathcal{P}_\lambda:= \mathcal{P}_{\lambda\kappa}$ denote a Poisson point
process of
intensity $\lambda\kappa$ on $[0,1]^d, d \geq2$, with $\kappa$ a bounded
density on $[0,1]^d$ and $\lambda\in(0, \infty)$. Given a closed subset
$\mathcal{M}\subset[0,1]^d$ of Hausdorff dimension $(d-1)$, we consider
general statistics $\sum_{x \in\mathcal{P}_\lambda} \xi(x, \mathcal{P}_\lambda, \mathcal{M})$, where the score
function $\xi$ vanishes unless the input $x$ is
close to $\mathcal{M}$ and where $\xi$ satisfies a weak spatial dependency
condition. We give a rate of normal convergence for the rescaled statistics
$\sum_{x \in\mathcal{P}_\lambda} \xi(\lambda^{1/d}x, \lambda^{1/d} \mathcal{P}_\lambda, \lambda
^{1/d} \mathcal{M})$ as $\lambda\to\infty$. When $\mathcal{M}$ is of class $C^2$, we
obtain weak laws of
large numbers and variance asymptotics for these statistics, showing
that growth
is surface order, that is, of order $\mathrm{Vol}(\lambda^{1/d} \mathcal{M})$. We use the
general results to deduce variance
asymptotics and central limit theorems for statistics arising in
stochastic geometry, including
Poisson--Voronoi volume and surface area estimators, answering
questions in
Heveling and Reitzner [\textit{Ann. Appl. Probab.} \textbf{19}
(2009) 719--736] and
Reitzner, Spodarev and Zaporozhets
[\textit{Adv. in Appl. Probab.} \textbf{44} (2012) 938--953].
The general results also yield the limit theory for the number of
maximal points in a sample.
\end{abstract}

%

\end{frontmatter}

\setcounter{footnote}{1}
\section{Main results}\label{INTRO}

\subsection{Introduction}\label{sec1}

Let $\P_\la:= \P_{\la\ka}$ denote a Poisson point process of
intensity $\la\ka$ on $[0,1]^d, d \geq2$, with $\ka$ a bounded
density on $[0,1]^d$ and $\la\in(0, \infty)$. Letting $\xi(\cdot,
\cdot)$ be a Borel measurable $\R$-valued function defined on pairs
$(x, \X)$, with $\X\subset\R^d$ finite and $x \in\X$, functionals
in stochastic
geometry may often be represented as linear statistics
$\sum_{x \in\P_\la} \xi(x, \P_\la)$. Here, $\xi(x, \P_\la)$
represents the contribution from~$x$,
which in general, depends on $\P_\la$. It is often more natural
to consider rescaled statistics
%
\begin{equation}
\label{ls} H^\xi(\P_\la):= \sum
_{x \in\P_\la} \xi \bigl(\la^{1/d}x, \la^{1/d}
\P_\la \bigr).
\end{equation}

Laws of large numbers, variance asymptotics and asymptotic
normality as $\la\to\infty$ for such statistics are established in
\cite{BY2,PeBer,PeEJP,PY1,PY5} with limits governed by the
behavior of $\xi$ at a point inserted
into the origin of a homogeneous Poisson point process. The sums $H^\xi
(\P_\la)$
exhibit growth of order $\Vol_d((\la^{1/d} [0,1])^d)= \la$, the
$d$-dimensional volume measure of the set carrying the scaled input
$\la^{1/d} \P_\la$. 
This gives the limit theory for score functions of nearest neighbor distances,
Voronoi tessellations, percolation and germ grain
models \cite{BY2,PeBer,PY1}. 
Problems of interest sometimes involve $\R$-valued score functions
$\xi$ of
three arguments, with the third being a set $\M\subset\R^d$ of
Hausdorff dimension $(d-1)$, and where scores $\xi(\la^{1/d}x, \la
^{1/d} \P_\la,
\la^{1/d}\M)$ vanish unless $x$ is close to $\M$.
This gives rise to
%
\begin{equation}
\label{lsM} H^\xi(\P_\la, \M):= \sum
_{x \in\P_\la} \xi \bigl(\la^{1/d}x, \la^{1/d}
\P_\la, \la ^{1/d}\M \bigr).
\end{equation}

Here,
$\M$ might represent the boundary of the support of $\ka$ or more
generally, the boundary of a fixed body, as would be the case in
volume and surface integral estimators. We show that modifications
of the methods used to study \eqref{ls} yield the limit theory of
\eqref{lsM}, showing that the scaling is surface order, that is,
$H^\xi(\P_\la, \M)$ is
order $\Vol_{d-1}(\la^{1/d}( \M\cap[0,1]^d))= \Theta(\la
^{(d-1)/d})$. The general limit theory for~\eqref{lsM}, as given in Section~\ref{sec1.2}, yields variance
asymptotics and central limit theorems
for the Poisson--Voronoi volume estimator, answering questions posed in
\cite{HR,RSZ}. 
We introduce a surface area estimator induced by Poisson--Voronoi tessellations
and we use the general theory to obtain its consistency and variance
asymptotics.
Finally, the general theory yields
the limit
theory for the number of maximal points in random sample, including
variance asymptotics and rates of normal
convergence, extending \cite{BHLT}--\cite{Ba}. See Section~\ref{applic} for
details. We anticipate further applications to germ-grain and
continuum percolation models, but postpone treatment of this.

\subsection{General results} \label{sec1.2}

We first introduce terminology, cf. \cite{BY2,PeBer,PeEJP,PY1,PY5}. Let ${\MM}(d)$ denote the collection of closed sets $\M\subset
[0,1]^d$ having finite $(d-1)$-dimensional Hausdorff measure.
Elements of ${\MM}(d)$ may or may not have
boundary and are endowed with the subset topology of $\R^d$. Let ${\MM
}_2(d) \subset{\MM}(d)$ denote those $\M\in
{\MM}(d)$ which are $C^2$, orientable submanifolds. Given $\M\in
{\MM}(d)$, almost all points $x \in
[0,1]^d$ are uniquely represented as
%
\begin{equation}
\label{param} x:= y + t {\mathbf{u}}_y,
\end{equation}
where $y:=y_x \in\M$ is the closest point in $\M$
to $x$, $t:=t_x \in\R$ and ${\mathbf{u}}_y$ is a fixed direction (see,
e.g., Theorem 1G of \cite{Fr}, \cite{HLW}); ${\mathbf{u}}_y$ coincides
with the unit outward normal to $\M$ at $y$ when $\M\in
{\MM}_2(d)$.
We write $x = (y_x,t_x):=(y,t)$ and shorthand $(y,0)$
as $y$ when the context is clear. To avoid pathologies, we assume
$\H^{d-1}(\M\cap
\partial([0,1]^d) ) =0$. Here, $\H^{d-1}$ denotes $(d-1)$-dimensional
Hausdorff measure,
normalized to coincide with $\Vol_{d-1}$ on hyperplanes.

Let $\xi(x, \X, \M)$ be a Borel measurable $\R$-valued function
defined on triples $(x, \X, \M)$, where $\X\subset\R^d$ is finite,
$x \in\X$, and $\M\in\MM(d)$. If $x \notin\X$, we shorthand
$\xi(x, \X\cup\{x\}, \M)$ as $\xi(x, \X, \M)$.
Let $S:= S(\M) \subset[0,1]^d$ be the set of points admitting the
unique representation
\eqref{param} and put $S':= \{(y_x, t_x)\}_{x \in S}$. If $(y,t) \in
S'$, then we put
$\xi((y,t), \X, \M) = \xi(x, \X, \M)$ where $x= y + t
{\mathbf{u}}_y$, otherwise we put $\xi((y,t), \X, \M) = 0$.

We assume $\xi$ is
translation invariant, that is, for all $z \in\R^d$ and input $(x, \X,
\M)$ we have
$ \xi(x, \X, \M) = \xi(x + z, \X+ z, \M+ z)$. Given $\la\in[1,
\infty)$, define dilated scores $\xi_\la$ by
%
\begin{equation}
\label{rescale} \xi_\la(x, \X, \M):= \xi \bigl(\la^{1/d}x,
\la^{1/d} \X, \la^{1/d} \M \bigr),
\end{equation}
so that
\eqref{lsM} becomes
%
\begin{equation}
\label{lsMn} H^\xi(\P_\la, \M):= \sum
_{x \in\P_\la} \xi_\la(x, \P_\la, \M).
\end{equation}
%

We recall two weak spatial dependence conditions for $\xi$.
For $\tau\in(0, \infty)$,
$\H_\tau$ denotes the homogeneous Poisson point process of intensity
$\tau$ on $\R^d$.
For all $x \in\R^d$, $r \in(0, \infty)$, let $B_r(x):= \{w \in\R
^d\dvtx \Vert x-w\Vert  \leq r
\}$, where $\Vert \cdot\Vert $ denotes Euclidean norm. Let $\0$ denote a point
at the origin of $\R^d$.
Say that $\xi$ is \emph{homogeneously stabilizing} if for all $\tau\in
(0, \infty)$
and all $(d-1)$-dimensional hyperplanes $\HH$, there is
$R:=R^\xi(\H_\tau,\HH) \in(0, \infty)$ a.s.
(a radius of stabilization) such
that
%
\begin{equation}
\label{hom} \xi \bigl(\0, \H_\tau\cap B_R(\0), \HH
\bigr) = \xi \bigl(\0, \bigl(\H_\tau\cap B_R(\0 ) \bigr)
\cup\A, \HH \bigr)
\end{equation}
for all locally finite $\A\subset B_R(\0)^c$.
Given \eqref{hom}, the definition of $\xi$ extends to infinite
Poisson input, that is, $\xi(\0, \H_\tau, \HH) = \lim_{r \to
\infty}
\xi(\0, \H_\tau\cap B_r(\0), \HH)$.

Given $\M\in\MM(d)$, say that
$\xi$ is \emph{exponentially stabilizing} with respect to the pair $(\P
_\la, \M)$ if
for all $x \in\R^d$ there is a radius of stabilization $R:=R^\xi(x,
\P_\la, \M) \in(0, \infty)$
a.s. such that
%
\begin{equation}
\label{expo} \xi_\la \bigl(x, \P_\la\cap
B_{\la^{-1/d}R}(x), \M \bigr) = \xi_\la \bigl(x, \bigl(\P
_\la\cap B_{\la^{-1/d}R}(x) \bigr) \cup\A, \M \bigr)
\end{equation}
for all locally finite $\A\subset\R^d \setminus B_{\la
^{-1/d}R}(x)$, and the tail probability
$\tau(t):= \tau(t, \M):= \sup_{\la> 0, x \in\R^d} P[R(x, \P_\la, \M) > t]$
satisfies $\limsup_{t \to\infty}
t^{-1} \log\tau(t) < 0$. 


Surface order
growth for the sums at \eqref{lsMn} should involve finiteness of the
integrated score $\xi_\la((y,t), \P_\la, \M)$ over $t \in\R$.
Thus, it is natural to require the following condition.
Given $\M\in\MM(d)$ and $p \in[1, \infty)$, say that $\xi$ \emph{satisfies the $p$ moment condition with respect to $\M$}
if there is a bounded integrable function $G^{\xi,p}:= G^{\xi, p, \M}\dvtx \R\to\R^+$ such that for all $u \in\R$
%
\begin{equation}
\label{mom} \sup_{z \in\R^d \cup\varnothing} \sup_{y \in\M} \sup
_{\la> 0} \E\bigl|\xi_\la \bigl( \bigl(y,
\la^{-1/d}u \bigr), \P_\la\cup z, \M \bigr)\bigr |^p
\leq G^{\xi,p}\bigl(|u|\bigr).
\end{equation}
%
Say that $\xi$ decays \emph{exponentially fast} with respect to the
distance to $\M$ if for all $p \in[1, \infty)$
%
\begin{equation}
\label{strongmom} \limsup_{|u| \to\infty} |u|^{-1} \log
G^{\xi,p}\bigl(|u|\bigr) < 0.
\end{equation}
%

Next, given $\M\in\MM_2(d)$ and $y \in\M$, let $\HH(y,\M)$ be
the $(d-1)$-dimensional hyperplane tangent to
$\M$ at $y$. 
Put $\HH_y:=\HH(\0,\M-y)$. The
score $\xi$ is \emph{well-approximated by $\P_\la$ input on
half-spaces} if for all $\M\in\MM_2(d)$, all $y \in\M$, and
all $w \in\R^d$, we have
%
\begin{eqnarray}
\label{lin}
&&\lim_{\la\to\infty} \E\bigl| \xi \bigl(w,
\la^{1/d}( \P_\la- y), \la^{1/d}(\M- y) \bigr)
\nonumber
\\[-8pt]
\\[-8pt]
\nonumber
&&\hspace*{66pt}{}- \xi
\bigl(w, \la^{1/d}(\P_\la- y), \HH_y \bigr)\bigr| =
0.
\end{eqnarray}
%

We now give three general limit theorems, proved in Sections~\ref{Proofs1} and \ref{sec5}. In Section~\ref{applic}, we use these
results to deduce the limit theory for statistics arising in
stochastic geometry. Let ${\cal C}(\M)$ denote the set of functions on
$[0,1]^d$ which are continuous at all points $y \in\M$.
Let $\0_y$ be a point at the origin of $\HH_y$.

\begin{theo}[(Weak law of large numbers)]\label{WLLN} Assume $\M\in
\MM_2(d)$
and $\ka\in{\cal C}(\M)$. If $\xi$ is homogeneously stabilizing
\eqref{hom}, satisfies
the moment condition \eqref{mom} for some $p > 1$, and is
well-approximated by
$\P_\la$ input on half-spaces \eqref{lin}, then
%
\begin{eqnarray}
\label{dWLLN} &&\lim_{\la\to\infty} \la^{-(d-1)/d}
H^\xi( \P_\la, \M)\nonumber\\
&&\qquad = \mu (\xi, \M)
\\
&&\hspace*{-3pt}\qquad:= \int_{\M} \int_{-\infty}^{\infty}
\E\xi \bigl((\0_y,u), \H_{\ka(y)}, \HH_y \bigr)
\ka(y) \,du \,dy \qquad\mbox{in } L^p.\nonumber
\end{eqnarray}
\end{theo}

Next, for $x, x' \in\R^d$, $\tau\in(0, \infty)$, and all
$(d-1)$-dimensional hyperplanes $\HH$ we put
\begin{eqnarray*}
&& c^\xi \bigl(x,x'; \H_\tau, \HH \bigr)\\
&&\qquad:= \E
\xi \bigl(x, \H_\tau\cup x', \HH \bigr) \xi
\bigl(x', \H_\tau\cup x, \HH \bigr) - \E\xi(x,
\H_\tau, \HH) \E\xi \bigl(x', \H_\tau, \HH
\bigr).
\end{eqnarray*}
Put for all $\M\in\MM_2(d)$
%
\begin{eqnarray}
\label{sigma} \sigma^2(\xi, \M)&:=& \mu \bigl(\xi^2, \M
\bigr)
\nonumber\\
&&{}+ \int_{\M} \int_{\R^{d-1}} \int
_{-\infty}^{\infty}\int_{-\infty}^{\infty}
c^\xi \bigl((\0_y,u), (z,s); \H_{\ka(y)},
\HH_y \bigr)\\
&&\hspace*{98pt}{}\times \ka(y)^2 \,du \,ds\,dz \,dy.\nonumber
\end{eqnarray}

\begin{theo}[(Variance asymptotics)] \label{Var} Assume $\M\in\MM
_2(d)$ and $\ka\in{\cal C}(\M)$. If~$\xi$ is homogeneously
stabilizing \eqref{hom},
exponentially stabilizing \eqref{expo}, satisfies
the moment condition \eqref{mom} for some $p > 2$, and is
well-approximated by
$\P_\la$ input on half-spaces \eqref{lin}, then
%
\begin{equation}
\label{dWLLM} \lim_{\la\to\infty} \la^{-(d-1)/d} \Var
\bigl[H^\xi(\P_\la, \M) \bigr] = \sigma^2(\xi,
\M) \in[0, \infty).
\end{equation}
\end{theo}

Let $N(0, \sigma^2)$ denote a mean zero normal random variable with
variance $\sigma^2$ and let $\Phi(t):= P[N(0,1) \leq t], t \in
\R$, be the distribution function
of the standard normal. 

\begin{theo}[(Rate of convergence to the normal)] \label{CLT}
Assume $\M\in\MM(d)$. If $\xi$ is exponentially stabilizing \eqref
{expo} and satisfies
exponential decay \eqref{strongmom} for some $p > q$, $q \in(2,3]$,
then there is a finite constant
$c:=c(d, \xi, p, q)$ such that for all $\la\geq2$
%
\begin{eqnarray}
\label{dCLT} &&\sup_{t \in\R} \biggl\llvert P \biggl[
\frac{ H^\xi(\P_\la, \M) - \E
[H^\xi(\P_\la, \M)]} { \sqrt{ \Var[H^\xi(\P_\la, \M)] } } \leq t \biggr] - \Phi(t) \biggr\rrvert
\nonumber
\\[-8pt]
\\[-8pt]
\nonumber
&&\qquad\leq c (\log\la)^{dq + 1} \la^{(d-1)/d} \bigl(\Var
\bigl[H^\xi(\P_\la, \M ) \bigr] \bigr)^{-q/2}.
\end{eqnarray}
In particular, if $\sigma^2(\xi, \M) > 0$, then putting $q = 3$
yields a rate of convergence
$O((\log\la)^{3d + 1} \la^{-(d-1)/2d})$ to the normal distribution.
\end{theo}


\begin{rems*}
(i) (\emph{Simplification of limits}.) If $\xi(x, \X, \M)$
is invariant under rotations of
$(x, \X, \M)$, then the limit $\mu(\xi, \M)$ at \eqref{dWLLN}
simplifies to
%
\begin{equation}
\label{rotinv} {}\mu(\xi, \M):= \int_{\M} \int
_{-\infty}^{\infty
} \E\xi \bigl((\0,u), \H_{\ka(y)},
\R^{d-1} \bigr) \,du\, \ka(y) \,dy,
\end{equation}
where 
$(\0,u) \in\R^{d-1} \times\R$. The limit \eqref{sigma} simplifies to
%
\begin{eqnarray}
\label{rotinv-1} \sigma^2(\xi, \M)&:=& \mu \bigl(\xi^2, \M
\bigr)
\nonumber\\
&&{}+ \int_{\M} \int_{\R^{d-1}} \int
_{-\infty}^{\infty}\int_{-\infty}^{\infty}
c^\xi \bigl((\0,u), (z,s); \H_{\ka(y)}, \R^{d-1}
\bigr) \\
&&\hspace*{97pt}{}\times \ka(y)^2 \,du \,ds\,dz \,dy.\nonumber
\end{eqnarray}
If, in addition, $\xi$ is homogeneous of order $\gamma$ in the sense
that for all $a \in(0, \infty)$ we have
\[
\xi\bigl(ax, a \cal X, \R^{d-1}\bigr) = a^\gamma  \xi\bigl(x, \cal X, \R^{d-1}\bigr),
\]
then putting
%
\begin{equation}
\label{Defmu} \mu(\xi, d):= \int_{-\infty}^{\infty} \E\xi
\bigl((\0,u), \H_{1}, \R^{d} \bigr) \,du
\end{equation}
we get
that $\mu(\xi, \M)$ further simplifies to
%
\begin{equation}
\label{simplestm} \mu(\xi, \M):= \mu(\xi, d-1) \int_{\M}
\kappa(y)^{(d-\gamma - 1)/d}
 \,dy.
\end{equation}
Similarly, the variance limit $\sigma^2(\xi, \M)$ simplifies to
\begin{eqnarray*}
\sigma^2(\xi, \M)&:=& \mu \bigl(\xi^2, d-1 \bigr) \int
_{\M}
 \kappa(y)^{(d-\gamma - 1)/d}
 \,dy
\\
&&{}+ \int_{\R^{d-1}} \int_{-\infty}^{\infty}
\int_{-\infty}^{\infty} c^\xi \bigl((\0,u),
(z,s); \H_{1}, \R^{d-1} \bigr) \,du \,ds \,dz\\
&&{}\times \int
_{\M}
\kappa(y)^{(d- 2\gamma - 2)/d}
\,dy.
\end{eqnarray*}
If $\ka\equiv1$, then putting
%
\begin{equation}
\label{Defnu} \nu(\xi, d):= \int_{\R^{d}} \int
_{-\infty}^{\infty}\int_{-\infty}^{\infty}
c^\xi \bigl((\0,u), (z,s); \H_{1}, \R^{d}
\bigr) \,du \,ds \,dz
\end{equation}
we get that \eqref{dWLLN} and \eqref{dWLLM}, respectively, reduce to
%
\begin{equation}
\label{supersimplemean} \lim_{\la\to\infty} \la^{-(d-1)/d}
H^\xi( \P_\la, \M) = \mu (\xi, d-1) \H^{d-1}(\M)
\qquad\mbox{in } L^p
\end{equation}
and
%
\begin{eqnarray}
\label{supersimplevar}&& \lim_{\la\to\infty} \la^{-(d-1)/d} \Var
\bigl[H^\xi(\P_\la, \M) \bigr]
\nonumber
\\[-8pt]
\\[-8pt]
\nonumber
&&\qquad = \bigl[\mu \bigl(
\xi^2, d -1 \bigr) + \nu(\xi, d-1) \bigr] \H^{d-1}(\M).
\end{eqnarray}

(ii) (\emph{A scalar central limit theorem}.)
Under the hypotheses of Theorems~\ref{Var} and~\ref{CLT}, we obtain
as $\la\to\infty$,
%
\begin{equation}
\label{scalarCLT} \la^{-(d-1)/2d} \bigl( H^\xi(\P_\la,\M)
- \E H^\xi(\P_\la,\M) \bigr) \stackrel{{\cal D}} {
\longrightarrow}N \bigl(0, \sigma^2(\xi, \M) \bigr).
\end{equation}
In general, separate arguments are needed to show strict positivity of\break
$\sigma^2(\xi, \M)$.

(iii) (\emph{Extensions to binomial input}.)
By coupling $\P_\la$ and binomial input $\{X_i\}_{i=1}^n$, where
$X_i, i \geq1$, are i.i.d.
with density $\ka$, it may be shown that
Theorems~\ref{WLLN} and~\ref{Var} hold for input $\{X_i\}_{i=1}^n$
under additional
assumptions on $\xi$. See Lemma \ref{dePoiss}.

(iv) (\emph{Extensions to random measures}.) Consider the
random measure
\[
\mu^\xi_\la:= \sum_{x \in\P_\la}
\xi_\la(x, \P_\la, \M) \delta_x,
\]
where $\delta_x$ denotes the Dirac point mass at $x$. For $f \in
\B([0,1]^d)$, the class of bounded functions on $[0,1]^d$, we put
$\langle f, \mu^\xi_\la\rangle:= \int f \,d \mu^\xi_\la$.
Modifications of the proof of Theorem \ref{WLLN} show that when $f
\in{\cal C}([0,1]^d)$, we have $L^p$, $p \in\{1, 2\}$, convergence
%
\begin{eqnarray}
\label{measLLN} &&\lim_{\la\to\infty} \la^{-(d-1)/d} \bigl\langle f,
\mu^\xi_\la \bigr\rangle\nonumber\\
&&\qquad= \mu(\xi,\M,f)
\\
&&\hspace*{-3pt}\qquad:= \int_{\M} \int_{-\infty}^{\infty}
\E\xi \bigl((\0_y,u), \H_{\ka
(y)}, \HH_y \bigr)
\ka(y) f(y) \,du \,dy.\nonumber
\end{eqnarray}
Using that a.e. $x \in[0,1]^d$ is a Lebesgue point for $f$, it may
be shown this limit extends to $f \in\B([0,1]^d)$ (Lemma 3.5 of
\cite{PeBer} and Lemma 3.5 of \cite{PeEJP}). The limit~\eqref{measLLN} shows up in surface integral approximation as seen
in Theorem \ref{SA} in Section~\ref{surface}.

Likewise, under the assumptions of Theorem \ref{Var}, it may be shown
for all $f \in\B([0,1]^d)$ that
\[
\lim_{\la\to\infty} \la^{-(d-1)/d} \Var \bigl[ \bigl\langle f,
\mu^\xi_\la \bigr\rangle \bigr] = \sigma^2(\xi,
\M,f),
\]
where
\begin{eqnarray*}
\sigma^2(\xi, \M,f)&:=& \mu \bigl(\xi^2,
\M,f^2 \bigr)
\\
&&{}+ \int_{\M} \int_{\R^{d-1}} \int
_{-\infty}^{\infty}\int_{-\infty}^{\infty}
c^\xi \bigl((\0_y,u), (z,s); \H_{\ka(y)},
\HH_y \bigr)\\
&&\hspace*{97pt}{}\times \ka(y)^2f(y)^2 \,du \,ds\,dz \,dy.
\end{eqnarray*}
Finally, under the assumptions of Theorem \ref{CLT}, we get the rate
of convergence
\eqref{dCLT} with $H^\xi(\P_\la, \M)$ replaced by $\langle f, \mu
^\xi_\la\rangle$.


(v) (\emph{Comparison with \cite{PY5}}.) Theorem
\ref{CLT} is the surface order analog of Theorem~2.1 of \cite{PY5}.
Were one to
directly apply the latter result to $H^\xi(\P_\la, \M)$, one would
get
\begin{eqnarray}
\label{uninform} &&\sup_{t \in\R} \biggl\llvert P \biggl[
{ H^\xi(\P_\la, \M) - \E H^\xi
(\P_\la, \M) \over\sqrt{ \Var
[H^\xi(\P_\la, \M)] } } \leq t \biggr] - \Phi(t) \biggr\rrvert
\nonumber
\\[-8pt]
\\[-8pt]
\nonumber
&&\qquad= O \bigl( (\log\la)^{3d + 1} \la \bigl( \Var \bigl[H^\xi(
\P_\la, \M) \bigr] \bigr)^{-3/2} \bigr).
\end{eqnarray}
However, when $\Var
[H^\xi(\P_\la, \M)] = \Omega(\la^{(d-1)/d})$, as is the case in
Theorem \ref{Var},
the right-hand side
of \eqref{uninform} is $O( (\log\la)^{3d + 1} \la^{-(d + 1)/2d }
)$. 
The reason for this suboptimal rate is that \cite{PY5}
considers sums of $\Theta(\la)$ nonnegligible contributions $\xi(x,
\P_\la)$, whereas here,
due to condition \eqref{strongmom}, the number of nonnegligible
contributions is surface order,
that is, of order $\Theta(\la^{(d-1)/d})$.

(vi) (\emph{Comparison with \cite{PY6}}.) Let $\M\in
\MM_2(d)$. In contrast with the present paper, \cite{PY6} considers
statistics $H^\xi(\Y_n):= \sum_{i=1}^n \xi(n^{1/(d-1)}Y_i,
n^{1/(d-1)} \Y_n)$, with input
$\Y_n:= \{Y_j \}_{j=1}^n$ carried by $\M$ rather than $[0,1]^d$.
In this set-up, $H^\xi(\Y_n)$ exhibits growth $\Theta(n)$.
\end{rems*}

\section{Applications} \label{applic}

\subsection{Poisson--Voronoi volume estimators}

Given $\P_\la$ as in Section~\ref{INTRO} and an unknown Borel set $A
\subset[0,1]^d$, suppose one can
determine which points in the realization of $\P_\la$ belong to $A$
and which belong to
$A^c:= [0,1]^d \setminus A$. How can one use this information to
establish consistent
statistical estimators of geometric properties of
$A$, including $\Vol(A)$ and $\H^{d-1}(\partial A)$? 
Here and henceforth, we shorthand $\Vol_d$ by $\Vol$.
In this section, we use our general results to give the limit theory
for a well-known estimator of $\Vol(A)$;
the next section proposes a new estimator of $\H^{d-1}(\partial A)$
and gives its limit theory as well.

For $\X\subset\R^d$ locally finite and $x \in\X$, let $C(x,
\X)$ denote the Voronoi cell generated by $\X$ and with center $x$.
Given $\P_\la$ 
and a Borel set $A \subset
[0,1]^d$, the Poisson--Voronoi
approximation of $A$  is the union of Voronoi cells with
centers inside $A$, namely
\[
{ A_\la:= \bigcup_{x \in\P_\la\cap A} C(x,
\P_\la).}
\]
The set $A_\la$ was introduced by Khmaladze and Toronjadze \cite{KT},
who anticipated that $A_\la$ should well-approximate the target $A$ in
the sense
that a.s. $\lim_{\la\to\infty} \Vol(A \Delta A_\la) = 0$.
This conjectured limit holds; as shown by \cite{KT} when $d = 1$
and by Penrose \cite{PeBer} for all $d = 1,2,\ldots.$
Additionally, if $\P_\la$ is replaced by a homogeneous Poisson point
process on $\R^d$ of intensity $\la$, then $\Vol(A_\la)$ is an
unbiased estimator of
$\Vol(A)$ (cf. \cite{RSZ}), 
rendering $A_\la$ of interest in image analysis, nonparametric
statistics and quantization,
as discussed in Section~1 of \cite{KT} as well as Section~1 of
Heveling and Reitzner \cite{HR}.

Heuristically, $\Vol(A_\la) - \E\Vol(A_\la)$ involves cell
volumes $\Vol(C(x, \P_\la))$, $x \in\P_\la$, within $O(\la
^{-1/d})$ of
$\partial A$. The number of such terms is of surface order, that is there
are roughly $O(\la^{(d-1)/d})$ such terms, each contributing
roughly $O(\la^{-2})$ toward the total variance. Were the terms spatially
independent, one might expect that as $\la\to\infty$,
%
\begin{equation}
\label{starr} \la^{(d + 1)/2d} \bigl(\Vol(A_\la) - \E
\Vol(A_\la) \bigr) \stackrel{{\cal D}} {\longrightarrow}N \bigl(0,
\sigma ^2 \bigr),
\end{equation}
%
as conjectured in Remark 2.2 of \cite{RSZ}.
We use Theorems \ref{Var}--\ref{CLT} to
prove this conjecture and to obtain a closed form expression for
$\sigma^2$ when $\partial A \in\MM_2(d)$; we find rates of normal
convergence for
$(\Vol(A_\la) - \E\Vol(A_\la))/ \sqrt{ \Var
\Vol(A_\la)}$ assuming only $\partial A \in\MM(d)$.
This adds to Schulte
\cite{Sch}, who for $\ka\equiv1$ and $A$ compact, convex, shows
that $ (\Var
\Vol(A_\la))^{-1/2} (\Vol(A_\la) - \E\Vol(A_\la)) $
is asymptotically normal, $\la\to\infty$. 
We obtain analogous limits for $\Vol(A \Delta A_\la)$.
In addition to the standing assumption $\Vert \ka\Vert _{\infty} < \infty$, we
assume everywhere in this section
that $\ka$ is bounded away from zero on $[0,1]^d$.


\begin{theo} \label{main1} If $\partial A \in\MM(d)$, then
\begin{eqnarray*}
&&\sup_{t \in\R} \biggl\llvert P \biggl[{\Vol(A_\la) - \E\Vol(A_\la)
\over\sqrt{ \Var
\Vol(A_\la)}}
\leq t \biggr] - \Phi(t) \biggr\rrvert \\
&&\qquad= O \bigl( (\log\la)^{3d + 1}
\la^{-2 - 1/d } \bigl( \Var \Vol( A_\la) \bigr)^{-3/2}
\bigr)
\end{eqnarray*}
and
\begin{eqnarray*}
&&\sup_{t \in\R} \biggl\llvert P \biggl[{\Vol(A \Delta A_\la) - \E\Vol(A
\Delta A_\la) \over\sqrt{ \Var
\Vol(A \Delta A_\la)}}
\leq t \biggr] - \Phi(t) \biggr\rrvert \\
&&\qquad= O \bigl( (\log\la)^{3d + 1}
\la^{-2 - 1/d } \bigl( \Var \Vol(A \Delta A_\la)
\bigr)^{-3/2} \bigr).
\end{eqnarray*}
\end{theo}

The rate of convergence is uninformative without lower bounds\break on $\Var
\Vol( A_\la)$ and $\Var
\Vol(A \Delta A_\la)$. Schulte \cite{Sch} shows $\Var\Vol( A_\la)
=\break \Omega(\la^{-(d + 1)/d})$ when $A$ is compact
and convex. The next result provides lower bounds when
$\partial A$ contains a smooth subset. 
For locally finite $\X\subset\R^d$, $x \in\X$, define the volume scores
%
\begin{eqnarray}
\label{defp}&& \nu^{ \pm}(x, \X, \partial A)
\nonumber
\\[-8pt]
\\[-8pt]
\nonumber
&&\qquad:= \cases{ %
\Vol \bigl(C(x, \X) \cap A^c \bigr), & \quad $\mbox{if }
C(x, \X) \cap\partial A \neq\varnothing, x \in A,$
\vspace*{2pt}\cr
\pm\Vol \bigl(C(x, \X) \cap A \bigr), & \quad
$\mbox{if } C(x, \X) \cap\partial A \neq
\varnothing, x \in A^c,$
\vspace*{2pt}\cr
0, & \quad$\mbox{if } C(x, \X) \cap\partial A = \varnothing.$}
\end{eqnarray}
In view of limits such as \eqref{rotinv-1}, we need to define scores
on hyperplanes $\R^{d-1}$. We thus put
%
\begin{eqnarray}\quad
\label{defpinf}&& \nu^{\pm} \bigl(x, \X, \R^{d-1} \bigr)
\nonumber
\\[-8pt]
\\[-8pt]
\nonumber
&&\qquad:=
\cases{ %
\Vol \bigl(C(x, \X) \cap\R_+^{d-1}
\bigr), &\quad  $\mbox{if } C(x, \X) \cap\R^{d-1} \neq\varnothing, x \in
\R_{-}^{d-1}$,
\vspace*{2pt}\cr
\pm\Vol \bigl(C(x, \X) \cap\R_{-}^{d-1} \bigr), &\quad $\mbox{if }
C(x, \X) \cap \R^{d-1} \neq\varnothing, x \in\R_+^{d-1}$,
\vspace*{2pt}\cr
0, & \quad $\mbox{if } C(x, \X) \cap\R^{d-1} = \varnothing,$}
\end{eqnarray}
where $\R_+^{d-1}:= \R^{d-1} \times[0, \infty)$ and
$\R_{-}^{d-1}:= \R^{d-1} \times(-\infty, 0]$. Define $\sigma^2(\nu
^-,\partial A)$ by putting
$\xi$ and $\M$ to be $\nu^-$ and $\partial A$, respectively, in
\eqref{rotinv-1}. Similarly, define
$\sigma^2(\nu^+,\partial A)$. When $\ka\equiv1$, these expressions
further simplify as at \eqref{supersimplevar}.

\begin{theo} \label{main2} If $\ka\in{\cal C}(\partial A)$
and if $\partial A$ contains a $C^1$ open subset,
then
\[
\Var \Vol( A_\la) = \Omega \bigl(\la^{-(d + 1)/d} \bigr)
\quad\mbox{and}\quad \Var \Vol(A \Delta A_\la) = \Omega \bigl(\la^{-(d + 1)/d}
\bigr).
\]
Additionally, if $\partial A \in\MM_2(d)$, then
\begin{eqnarray*}
\lim_{\la\to\infty} \la^{(d + 1)/d} \Var\Vol(A_\la) &=&
\sigma ^2 \bigl(\nu^-,\partial A \bigr) \quad\mbox{and}\\
 \lim
_{\la\to\infty} \la^{(d + 1)/d} \Var\Vol(A \Delta A_\la)
&=& \sigma^2 \bigl(\nu^+,\partial A \bigr).
\end{eqnarray*}
%
\end{theo}

Combining the above results gives the following central limit theorem
for $\Vol(A_\la) - \Vol(A)$;
identical results hold for $\Vol(A \Delta A_\la) - \E\Vol(A \Delta
A_\la)$.

\begin{coro} \label{main3} If $\ka\in{\cal C}(\partial A)$ and
if either $\partial A$ contains a $C^1$ open subset or $A$ is
compact and convex, then
\[
\sup_{t \in\R} \biggl\llvert P \biggl[{\Vol(A_\la) - \E\Vol(A_\la
)\over\sqrt{ \Var
\Vol(A_\la)}}
\leq t \biggr] - \Phi(t) \biggr\rrvert = O \bigl( (\log\la )^{3d + 1}
\la^{-(d - 1)/2d } \bigr).
\]
Additionally, if $\partial A \in\MM_2(d)$, then as $\la\to\infty$
\[
\la^{(d + 1)/2d} \bigl( \Vol(A_\la) - \E\Vol(A_\la)
\bigr) \stackrel{{\cal D}} {\longrightarrow}N \bigl(0, \sigma ^2
\bigl( \nu^-, \partial A \bigr) \bigr).
\]
\end{coro}

Recall $X_i, i \geq1$, are i.i.d. with density $\ka$; $\X_n:= \{X_i\}
_{i = 1}^n$. The binomial-Voronoi approximation of $A$ is $
A_n:= \bigcup_{X_i \in A} C(X_i, \X_n)$. The above theorems extend
to binomial input as follows.


%
\begin{theo} \label{main4ab} If $\ka\in{\cal C}(\partial A)$ and
if either $\partial A$ contains a $C^1$ open subset or $A$ is
compact and convex, then
\[
\Var \Vol( A_n) = \Omega \bigl(n^{-(d + 1)/d} \bigr)\quad \mbox{and}\quad
\Var \Vol(A \Delta A_n) = \Omega \bigl(n^{-(d + 1)/d} \bigr).
\]
Additionally, if $\partial A \in\MM_2(d)$, then
\begin{eqnarray*}
\lim_{n \to\infty} n^{(d + 1)/d} \Var\Vol(A_n) &=&
\sigma^2 \bigl(\nu ^-,\partial A \bigr), \\
\lim_{n \to\infty}
n^{(d + 1)/d} \Var\Vol(A \Delta A_n) &=& \sigma^2 \bigl(
\nu^+,\partial A \bigr),
\end{eqnarray*}
and as $n \to\infty$,
\[
n^{(d + 1)/2d} \bigl( \Vol(A_n) - \E\Vol(A_\la) \bigr)
\stackrel{{\cal D}} {\longrightarrow}N \bigl(0, \sigma^2 \bigl(\nu
^-, \partial A \bigr) \bigr).
\]
\end{theo}


\begin{rems*}
(i) (\emph{Theorem \ref{main2}}.) When $\ka\equiv1$,
Theorem \ref{main2} and \eqref{supersimplevar} show that the limiting
variance of
$\Vol( A_\la)$ and $\Vol(A \Delta A_\la)$ involve multiples of $\H
^{d-1}(\partial A)$,
settling a conjecture implicit in Remark 2.2 of \cite{RSZ} when
$\partial A \in\MM_2(d)$.
Up to now, it has been known that $\Var\Vol(A_\la) = \Theta(\la
^{-(d +
1)/d})$ for $A$ compact and convex, where the upper and lower bounds
follow from
\cite{HR} and \cite{Sch}, respectively.

(ii) (\emph{Corollary \ref{main3}}.) When $\partial A$
contains a $C^1$ open subset, Corollary \ref{main3} answers the
first conjecture in Remark 2.2 of \cite{HR}; when $A$ is convex it
establishes a
rate of normal convergence for ${(\Vol(A_\la) - \E\Vol(A_\la)) /
\sqrt{
\Var\Vol(A_\la)}} $,
extending the main result of \cite{Sch} (Theorem 1.1). 

(iii) (\emph{The $C^2$ assumption}.) 
If $A \subset\R^d$ has finite perimeter, denoted $\operatorname{Per}(A)$, then~\cite{RSZ} shows that $\lim_{\la\to\infty} \la^{1/d} \E\Vol(A
\Delta A_\la)
= c_d \operatorname{Per}(A)$, where $c_d$ is an explicit constant depending only
on dimension. This remarkable result,
based on covariograms, holds with no other assumptions on $A$. Theorem
\ref{main2}
and Corollary~\ref{main3} hold for $\partial A$ not necessarily in
$\MM_2(d)$;
see \cite{TY}.
\end{rems*}

\subsection{Poisson--Voronoi surface integral estimators} \label{surface}

We show that the surface area of 
$A_\la$, when corrected by a factor independent of $A$,
consistently estimates $\H^{d-1}(\partial A)$ and that it satisfies
the limits in Theorems \ref{WLLN}--\ref{CLT}.

Given $\X$ locally finite and a Borel subset $A \subset\R^d$, define
for $x \in\X\cap A$
the area score $\alpha(x, \X, \partial A)$ to be the $\H^{d-1}$
measure of the $(d-1)$-dimensional faces of $C(x,\X)$ belonging to the
boundary of $\bigcup_{w \in\X\cap A} C(w, \X)$; if there are no
such faces or
if $x \notin\X\cap A$,
then set $\alpha(x, \X, \partial A)$ to be zero.
Similarly, for $x \in\X\cap\R^{d-1}_{-}$, put $\alpha(x, \X, \R
^{d-1})$ to be
the $\H^{d-1}$
measure of the $(d-1)$-dimensional faces of $C(x,\X)$ belonging to the
boundary of $\bigcup_{w \in\X\cap\R^{d-1}_{-}} C(w, \X)$,
otherwise $\alpha(x, \X, \R^{d-1})$ is zero.

The surface area of $A_\la$ is then given by $\sum_{x \in\P_\la}
\alpha(x, \P_\la, \partial A)$.
We might expect that the statistic
%
\begin{equation}
\label{alphastat} \la^{-(d-1)/d} H^\alpha(\P_\la, \partial
A) = \la^{-(d-1)/d}\sum_{x \in\P_\la}
\alpha_\la(x, \P_\la, \partial A)
\end{equation}
consistently estimates $\H^{d-1}(\partial A), \la\to\infty$, and more
generally, for $f \in B([0,1]^d)$ that
\[
\la^{-(d-1)/d} \sum_{x \in\P_\la} \alpha_\la(x,
\P_\la, \partial A) f(x)
\]
consistently estimates the surface integral $\int_{\partial
A} f(x) \H^{d-1}(dx)$. Provided that one introduces a \emph{universal
correction factor which is independent of the target $A$}, this
turns out to be the case, as seen in the next theorem. Define $\mu
(\alpha, d)$ and $\nu(\alpha, d)$
by putting $\xi$ to be $\alpha$ in
\eqref{Defmu} and \eqref{Defnu}, respectively.

\begin{theo} \label{SA} If $\ka\equiv1$ and $\partial A \in\MM
_2(d)$, then
%
\begin{equation}
\label{areaWLLN} \lim_{\la\to\infty} \bigl(\mu(\alpha, d-1)
\bigr)^{-1} \H^{d-1} (\partial A_\la) =
\H^{d-1} (\partial A) \qquad\mbox{in } L^2
\end{equation}
and
%
\begin{eqnarray}
\label{areaVar} &&\lim_{\la\to\infty} \la^{(d-1)/d} \Var \bigl[
\H^{d-1} (\partial A_\la ) \bigr]
\nonumber
\\[-8pt]
\\[-8pt]
\nonumber
&&\qquad= \bigl[\mu \bigl(
\alpha^2, d-1 \bigr) + \nu(\alpha, d - 1) \bigr]\H^{d-1} (
\partial A).
\end{eqnarray}
Further, for $f \in B([0,1]^d)$
%
\begin{eqnarray}
\label{areaWLLN1} &&\lim_{\la\to
\infty} \bigl(\mu(\alpha, d-1)
\bigr)^{-1} \la^{-(d-1)/d} \sum_{x \in\P_\la}
\alpha_\la(x, \P_\la, \partial A) f(x)
\nonumber
\\[-8pt]
\\[-8pt]
\nonumber
&&\qquad = \int
_{\partial A} f(x) \H^{d-1}(dx) \qquad\mbox{in } L^2.
\end{eqnarray}
\end{theo}

\begin{rems*}
(i) {(\textit{Extensions}.)} Assuming only $\partial A \in\MM
(d)$, it follows\break from Theorem \ref{CLT} and the upcoming proof of
Theorem \ref{SA} that\break $(\Var\H^{d-1}(\partial A_\la))^{-1/2}\times  (\H
^{d-1}(\partial A_\la) - \E\H^{d-1}(\partial
A_\la))$ is asymptotically normal.
When $\partial A \in\MM_2(d)$ it follows by
\eqref{scalarCLT} that as $\la\to\infty$
\[
\la^{-(d-1)/2d} \bigl(\H^{d-1}(\partial A_\la) - \E
\H^{d-1}(\partial A_\la) \bigr) \stackrel{{\cal D}} {
\longrightarrow}N \bigl(0, \sigma^2 \bigr),
\]
with
$\sigma^2:= [\mu(\alpha^2, d-1) + \nu(\alpha, d - 1)]\H^{d-1}
(\partial A)$.
Analogs of \eqref{areaWLLN}--\eqref{areaWLLN1} hold if
$\P_\la$ is replaced by $\X_n:= \{X_i\}_{i=1}^n$, $A_\la$ is
replaced by $A_n:= \bigcup_{X_i \in A} C(X_i, \X_n)$, and $n \to
\infty$.


(ii) {(\textit{Related work}.)} Using the Delaunay triangulation
of $\P_\la$,
\cite{JY} introduces an a.s. consistent estimator of surface
integrals of possibly nonsmooth boundaries. The limit theory for
the Poisson--Voronoi estimator $H^{\alpha}(\P_\la, \partial A)$
extends to non\-smooth $\partial A$ as in \cite{TY}.
\end{rems*}

\subsection{Maximal points} \label{maxsection}

Let $K \subset\R^d$ be a cone with nonempty interior and apex at
the origin of $\R^d$. Given $\X\subset\R^d$ locally finite, $x \in
\X$ is called $K$-maximal, or simply maximal if $(K \oplus x) \cap
\X= x$. Here, $K \oplus x$ is Minkowski addition, namely $K \oplus x:=
\{z + x\dvtx z \in K\}$. In the case $K = (\R^+)^d$, a point $x=
(x_1,\ldots,x_d) \in
\X$ is maximal if there is no other point $(z_1,\ldots,z_d) \in\X$
with $z_i \geq x_i$ for all $1 \leq i \leq d$. The maximal layer
$m_K(\X)$ is the  collection of maximal points in $\X$. Let
$M_K(\X):= \operatorname{card} (m_K(\X))$.

Maximal points feature in various disciplines. They are of broad
interest in computational geometry; see
books by Preparata and Shamos \cite{PS}, 
Chen et al.~\cite{Ch}.
Maximal points appear in pattern classification,
multicriteria decision analysis, networks, data mining, analysis of
linear programming and statistical decision theory; see
Ehrgott \cite{Eh} and Pomerol and Barba-Romero \cite{PB}.
In economics, when $K=(\R^+)^d$, the maximal layer and $K$ are termed
the Pareto set
and Pareto cone, respectively;
see Sholomov \cite{Sh} for a survey on Pareto optimality.

Next, let $\ka$ be a density having support
\[
A:= \bigl\{(v,w)\dvtx v \in D, 0 \leq w \leq F(v) \bigr\},
\]
where $F\dvtx D \to\R$ has continuous partials $F_i, 1 \leq i \leq d
-1$, which are bounded away from zero and negative infinity; $D \subset
[0,1]^{d-1}$, and $|F| \leq1$. Let $\P_\la:= \P_{\la\ka}$ and
$\X_n:= \{X_i\}_{i=1}^n$ as above.

Using Theorems \ref{WLLN}--\ref{CLT}, we deduce
laws of large numbers, variance asymptotics, and
central limit theorems for $M_K(\P_\la)$ and $M_K(\X_n)$, as $\la
\to\infty$ and $n \to\infty$, respectively. Put $\partial A:=
\{(v, F(v))\dvtx v \in D\}$ and let
\begin{eqnarray*}
\zeta(x, \X, \partial A):= \cases{ %
1, &\quad $\mbox{if
} \bigl((K \oplus x) \cap A \bigr) \cap\X= x,$
\vspace*{2pt}\cr
0, & \quad$\mbox{otherwise}.$}
\end{eqnarray*}
When $x = (y,t), y \in\partial A$, we write
%
\begin{equation}
\label{defz} \zeta(x, \X, \HH_y):= \cases{ %
1, & \quad$\mbox{if } \bigl((K \oplus x) \cap\HH_+(y, \partial A)
\bigr) \cap\X= x,$
\vspace*{2pt}\cr
0, &\quad $\mbox{otherwise},$}
\end{equation}
where $\HH_+(y, \partial A)$ is the half-space
containing $\0$ and with hyperplane $\HH(y, \partial A)$.

To simplify the presentation, we take $K = (\R^+)^d$, but the
results extend to general cones. Recalling definitions
\eqref{dWLLN} and \eqref{sigma}, we have the following results. 

\begin{theo} \label{main4} If $\ka\in{\cal C}(\partial A)$ and if
$\ka$ is bounded away from $0$ on $A$,
then
\begin{eqnarray}\label{zetaLLN}
&&\lim_{\la\to\infty} \la^{-(d-1)/d} M_K(
\P_\la)\nonumber\\
&&\qquad= \mu(\zeta, \partial A)
\\
 &&\qquad= ({d!})^{1/d} d^{-1} \Gamma
\bigl(d^{-1} \bigr) \int_D \Biggl\llvert \prod
_{i=1}^{d-1} F_i(v) \Biggr
\rrvert ^{1/d} \ka \bigl(v, F(v) \bigr)^{(d-1)/d} \,dv \qquad\mbox{in }
L^2\nonumber
\end{eqnarray}
and
%
\begin{equation}
\label{zetaVar} \lim_{\la\to\infty} \la^{-(d-1)/d} \Var
\bigl[M_K(\P_\la) \bigr] = \sigma^2(\zeta,
\partial A) \in(0, \infty).
\end{equation}
Moreover, as $\la\to\infty$, we have
\[
\la^{-(d-1)/2d} \bigl(M_K(\P_\la) - \E
M_K(\P_\la) \bigr) \stackrel{{\cal D}} {\longrightarrow}N
\bigl(0, \sigma ^2(\zeta, \partial A) \bigr).
\]
Identical limits hold with $M_K(\P_\la)$ replaced by $M_K(\X_n)$, $n
\to\infty$.
We also have
%
\begin{equation}
\label{maxCLT}\qquad \sup_{t \in\R} \biggl\llvert P \biggl[
\frac{ M_K(\P_\la) - \E M_K(\P
_\la)} { \sqrt{ \Var[M_K(\P_\la)] } } \leq t \biggr] - \Phi(t) \biggr\rrvert \leq c (\log
\la)^{3q + 1} \la^{(d-1)/2d}.
\end{equation}
\end{theo}

\begin{rems*}
(i) (\emph{Related expectation and
variance asymptotics}.) Formu\-la~\eqref{zetaVar} is new for all
dimensions $d$, whereas formula \eqref{zetaLLN} is new for $d > 2$.
For $d = 2$, \eqref{zetaLLN} extends work of Devroye \cite{De}, who
treats the case $\ka\equiv1$. Barbour and Xia \cite{BX,BX1}
establish growth rates for $ \Var[M_K(\P_\la)]$ but do not
determine limiting means or variances for $d > 2$. Hwang and Tsai \cite{HT}
determine $\E M_K(\X_n)$ and $\Var M_K(\X_n)$ when $A:= \{ (x_1,\ldots,x_d)\dvtx x_i \geq0,  \sum_{i=1}^d x_i \leq1 \}$, that is, $\partial A$ is
a subset
of the plane $\sum_{i=1}^d x_i = 1$.

(ii) (\emph{Related central limit theorems}.) Using Stein's
method, Barbour and Xia \cite{BX,BX1}
show for $d = 2$, $\ka$ uniform and $K = (\R^+)^2$ that
$(M_K(\X_n) - \E M_K(\X_n))/ \sqrt{ \Var M_K(\X_n)}$ tends to a
standard normal.
Assuming differentiability conditions on $F$, they find
rates of normal convergence of $M_K(\X_n)$ and $M_K(\P_\la)$ with
respect to the bounded Wasserstein distance \cite{BX} and the
Kolmogorov distance~\cite{BX1}, respectively.
Their work
adds to Bai et al. \cite{BHLT}, which for $K = (\R^+)^2$ establishes
variance asymptotics and central limit theorems when $\ka$ is
uniform on a convex polygonal region, and Baryshnikov \cite{Ba}, who
proves a central limit theorem under general conditions
on $\partial A$, still in the setting of homogeneous point sets. %

(iii) (\emph{Related results}.)
Parametrizing points
in $\R^d$ with respect to a fixed $(d-1)$-dimensional plane $\HH_0$,
the preprint \cite{BYpre} obtains expectation and variance asymptotics for
$M_K(\P_\la)$ and $M_K(\X_n)$, with limits depending on an integral
over the projection of $\partial A$ onto $\HH_0$. By comparison,
the limits in Theorem \ref{main4} follow straightforwardly from the
general limit theorems and exhibit an explicit dependence on the
graph of $F$, that is, $\partial A$. Preprint \cite{BYpre} uses
cumulants to show asymptotic normality without delivering the rate
of convergence offered by Theorem \ref{CLT}.

(iv) (\emph{Extensions}.) Separate analysis is needed to
extend Theorem \ref{main4} to spherical boundaries $\mathbb{S}^{d-1}
\cap[0, \infty)^d$, that is to say quarter circles in $d = 2$.
\end{rems*}

\subsection{Navigation in Poisson--Voronoi tessellations}

Put $\ka\equiv1$. Let $\X\subset\R^2$ be locally finite and let
$r(t), 0 \leq t \leq1$, be a $C^1$ curve
$\C$ in $[0,1]^2$. Let $\V_{\C}:= \V_{\C}(\X)$ be the union of
the Voronoi cells
$C(x, \X)$ meeting $\C$. Order the constituent cells of
$\V_{\C}$ according to the ``time'' at which $r(t)$ first meets the
cells. Enumerate the cells as
\[
C(x_1,\X,\C),\ldots,C(x_N,\X,\C); \qquad N \mbox{ random}.
\]
The piecewise linear path joining the nodes $x_1,\ldots,x_N$ is a
path $\C(\X)$ whose length $| \C(\X) |$ approximates the
length of $\C$. The random path $\C(\P_\la)$ has been studied by Bacelli
et al. \cite{BTZ}, which restricts to linear $\C$. For all $x \in
\X$ define the score
\[
\rho(x, \X, \C):= \cases{ %
\mbox{one half the sum
of lengths of edges incident to $x$ in}\vspace*{2pt}\cr
\qquad \C(\X) \mbox{ if } x \in\C(\X),
\vspace*{2pt}\cr
0,\qquad \mbox{otherwise}.}
\]
Then the path length $| \C(\P_\la) |$ satisfies
\[
\bigl| \C(\P_\la) \bigr| = \sum_{x \in\P_\la} \rho(x,
\P_\la, \C) = \la^{-1/2} H^{\rho}(\P_\la,
\C).
\]
We claim that the score $\rho$ satisfies the conditions of Theorems
\ref{WLLN}--\ref{CLT} and that therefore the limit theory of $| \C
(\P_\la)|$ may be
deduced from these general theorems, adding to \cite{BTZ}. Likewise,
using the Delaunay triangulation of $\P_\la$, one can find a unique
random path $\tilde{\C}_\la(\P_\la)$ whose edges meet $\C$ and
belong to the
triangulation of $\P_\la$, with length
\[
\bigl|\tilde{\C}_\la(\P_\la)\bigr| = \sum
_{x \in\P_\la} \tilde{\rho }(x, \P_\la, \C) =
\la^{-1/2} H^{\tilde{\rho}}(\P_\la, \C),
\]
where
\begin{eqnarray*}
&& \tilde{\rho}(x, \P_\la, \C)
\\
&&\qquad := \cases{ %
\mbox{one half the sum of lengths of edges incident to $x$ } \mbox{if } x \in
\tilde{\C}_\la(\P_\la),
\vspace*{2pt}\cr
0, \qquad\mbox{otherwise}.}
\end{eqnarray*}
Theorems \ref{WLLN}--\ref{CLT} provide the limit
theory for $|\tilde{\C}_\la(\P_\la)|$.
\section{Auxiliary results}

We give three lemmas pertaining to the rescaled scores
$\xi_\la, \la> 0$, defined at \eqref{rescale}.

\begin{lemm} \label{L1} Fix $\M\in\MM_2(d)$. Assume that $\xi$ is
homogeneously stabilizing, satisfies the
moment condition \eqref{mom} for $p > 1$ and is well-approximated by
$\P_\la$ input on
half-spaces \eqref{lin}. Then for almost all $y \in\M$, all $u \in
\R$, and all $x \in\R^d \cup\varnothing$
we have
%
\begin{equation}
\label{L3.1f}\qquad \lim_{\la\to\infty} \E\xi_\la \bigl(
\bigl(y, \la^{-1/d}u \bigr) + \la^{-1/d}x, \P_\la, \M
\bigr) = \E \xi \bigl((\0_y,u) + x, \H_{\ka(y)},
\HH_y \bigr).
\end{equation}
\end{lemm}

\begin{pf} Fix $\M\in\MM_2(d)$. We first show for
almost all $y \in\M$ that there exist coupled realizations $\P'_\la$
and $\H'_{\ka(y)}$ of $\P_\la$ and $\H'_{\ka(y)}$, respectively,
such that for $u \in\R$ and
$x \in\R^d$, we have as $\la\to\infty$
%
\begin{equation}
\label{conv1a} \xi_\la \bigl( \bigl(y, \la^{-1/d}u \bigr) +
\la^{-1/d}x, \P'_\la, \M \bigr) \stackrel {{\cal
D}} {\longrightarrow}\xi \bigl((\0_y,u) + x, \H'_{\ka(y)},
\HH_y \bigr).
\end{equation}
By translation invariance of $\xi$, we have
\begin{eqnarray*}
\xi_\la \bigl( \bigl(y, \la^{-1/d}u \bigr) +
\la^{-1/d}x, \P_\la, \M \bigr)& =& \xi_\la \bigl(
\bigl(\0_y, \la^{-1/d}u \bigr) + \la^{-1/d}x,
\P_\la- y, \M- y \bigr)
\\
&=& \xi \bigl((\0_y,u) + x, \la^{1/d}(\P_\la-
y), \la^{1/d}(\M- y) \bigr).
\end{eqnarray*}
By the half-space approximation assumption \eqref{lin},
we need only show for almost all $y \in\M$ that
there exist coupled realizations $\P'_\la$ and $\H'_{\ka(y)}$ of
$\P_\la$ and $\H_{\ka(y)}$, respectively,
such that as $\la\to\infty$
%
\begin{equation}
\label{conv1} \xi \bigl((\0_y,u) + x, \la^{1/d} \bigl(
\P'_\la- y \bigr), \HH_y \bigr) \stackrel{{
\cal D}} {\longrightarrow}\xi \bigl((\0_y,u) + x,
\H'_{\ka(y)}, \HH_y \bigr).
\end{equation}
This, however, follows from the homogeneous stabilization of $\xi$ and
the continuous mapping theorem; see Lemmas 3.2 and 3.2 of
\cite{PeBer}, which proves this assertion for the more involved case
of binomial input. 
Thus, \eqref{conv1a} holds and Lemma~\ref{L1} follows from uniform
integrability of $\xi_\la((y,
\la^{-1/d}u) + \la^{-1/d}x, \P'_\la, \M)$, which follows from the
moment condition \eqref{mom}.
\end{pf}


\begin{lemm} \label{L2} Fix $\M\in\MM_2(d)$. Assume that $\xi$ is
homogeneously stabilizing, satisfies the
moment condition \eqref{mom} for $p > 2$, and is well-approximated by
$\P_\la$ input on half-spaces \eqref{lin}.
Given $y \in\M$, $x \in\R^d$ and $u \in\R$,
put
\begin{eqnarray*}
X_\la&:=& \xi_\la \bigl( \bigl(y, \la^{-1/d}u
\bigr), \P_\la\cup \bigl( \bigl(y, \la^{-1/d}u \bigr) +
\la^{-1/d}x \bigr), \M \bigr),
\\
Y_\la&:=& \xi_\la \bigl( \bigl(y, \la^{-1/d}u
\bigr)+ \la^{-1/d}x, \P_\la\cup \bigl(y, \la^{-1/d}u
\bigr), \M \bigr),
\\
X&:=& \xi \bigl((\0_y,u), \H_{\ka(y)} \cup \bigl((
\0_y,u) + x \bigr), \HH_y \bigr) \quad\mbox{and}
\\
Y&:=& \xi \bigl((\0_y,u) + x, \H_{\ka(y)} \cup(
\0_y,u), \HH_y \bigr).
\end{eqnarray*}
Then for almost all $y \in\M$ we have $ \lim_{\la\to\infty} \E
X_\la Y_\la= \E X Y$.
\end{lemm}

\begin{pf} By the moment condition \eqref{mom}, the
sequence $X_\la^2, \la\geq1$, is uniformly integrable and hence
the convergence in distribution $X_\la\stackrel{{\cal
D}}{\longrightarrow}X$ extends to $L^2$
convergence and likewise for $Y_\la\stackrel{{\cal
D}}{\longrightarrow}Y$. The triangle
inequality and the Cauchy--Schwarz inequality give
\[
\Vert X_\la Y_\la- XY \Vert _1 \leq\Vert
Y_\la\Vert _2 \Vert X_\la- X \Vert
_2 + \Vert X\Vert _2 \Vert Y_\la- Y\Vert
_2.
\]
Lemma \ref{L2} follows since $\sup_{\la> 0} \Vert Y_\la\Vert _2 < \infty$
and $\Vert X\Vert _2 < \infty$.
\end{pf}

The next result quantifies the exponential decay of correlations
between scores on re-scaled input separated by Euclidean distance
$\Vert x\Vert $.

\begin{lemm} \label{L3} Fix $\M\in\MM(d)$. Let $\xi$ be
exponentially stabilizing \eqref{expo} and assume the moment condition
\eqref{mom} holds for some $p > 2$.
Then there is a $c_0 \in(0, \infty)$ such that for all $w, x \in\R^d$
and $\la\in(0, \infty)$, we have
\begin{eqnarray*}
&&\bigl|\E\xi_\la \bigl(w, \P_\la\cup \bigl(w +
\la^{-1/d}x \bigr), \M \bigr) \xi_\la \bigl(w +
\la^{-1/d}x, \P_\la\cup w , \M \bigr) \\
&&\hspace*{76pt}{}-
\E\xi_\la(w, \P_\la, \M) \E\xi_\la \bigl(w +
\la^{-1/d}x, \P_\la, \M \bigr)\bigr| \\
&&\qquad\leq c_0 \exp
\bigl(-c_0^{-1} \Vert x\Vert \bigr).
\end{eqnarray*}
\end{lemm}

\begin{pf} See the proof of Lemma 4.2 of \cite{PeEJP} or
Lemma 4.1 of \cite{BY2}.
\end{pf}

\section{Proofs of Theorems \texorpdfstring{\protect\ref{WLLN}--\protect\ref{Var}}{1.1--1.2}}
\label{Proofs1}

Roughly speaking, putting $x = \varnothing$ in \eqref{L3.1f} and integrating
\eqref{L3.1f} over $y \in\M$ and $u \in\R$, we obtain expectation
convergence of
$\la^{-(d-1)/d} H^\xi(\P_\la, \M)$ in Theorem \ref{WLLN}. We then
upgrade this to $L^1$ and $L^2$ convergence.
Regarding Theorem \ref{Var},
Lemmas \ref{L1} and \ref{L2} similarly yield convergence of the
covariance of scores
$\xi_\la$ at points $(y, \la^{-1/d}u)$ and $(y, \la^{-1/d}u) + \la^{-1/d}x$
and Lemma \ref{L3}, together with dominated convergence, imply
convergence of
integrated covariances over $x \in\R^d$ and $u \in\R$, as they
appear in the
iterated integral formula for $\la^{-(d-1)/d} \Var H^\xi(\P_\la, \M
)$. The details
go as follows.

\begin{pf*}{Proof of Theorem \ref{WLLN}} We first prove
$L^2$ convergence. Recall the definitions of $H^\xi(\P_\la,\M)$ and
$\mu(\xi, \M)$ at \eqref{lsM} and \eqref{dWLLN}, respectively. In
view of the identity
\begin{eqnarray*}
&&\E \bigl(\la^{-(d-1)/d } H^\xi(\P_\la,\M) - \mu(\xi,
\M) \bigr)^2
\\
&&\qquad= \la^{-2(d-1)/d} \E H^\xi(\P_\la,\M)^2 -
2 \mu(\xi, \M) \la^{-(d-1)/d} \E H^\xi(\P_\la,\M)\\
&&\qquad\quad{} +
\mu(\xi , \M)^2,
\end{eqnarray*}
it suffices to show
%
\begin{equation}
\label{one} \lim_{\la\to\infty} \la^{-(d-1)/d} \E
H^\xi(\P_\la,\M) = \mu(\xi, \M)
\end{equation}
and
%
\begin{equation}
\label{two} \lim_{\la\to\infty} \la^{-2(d-1)/d} \E
H^\xi(\P_\la,\M)^2 = \mu(\xi,
\M)^2.
\end{equation}
%
To show \eqref{one}, we first write
\[
\la^{-(d-1)/d} \E H^\xi(\P_\la,\M) =
\la^{1/d} \int_{[0,1]^d} \E \xi_\la(x,
\P_\la,\M) \ka(x)\,dx.
\]
Given $\M\in\MM_2(d)$ and $x \in[0,1]^d$, recall from \eqref{param}
the parameterization $x = y + t {\mathbf{u}}_y$, with ${\mathbf{u}}_y$ the
unit outward normal
to $\M$ at $y$. The Jacobian
of the map $h\dvtx x \mapsto(y + t {\mathbf{u}}_y)$ 
at $(y,t)$ is $J_h((y,t)):= \prod_{i=1}^{d-1}(1 +
t C_{y,i})$, where $C_{y,i}, 1 \leq i \leq d - 1$, are the principal
curvatures of $\M$ at $y$. Surfaces in $\MM_2(d)$ have
bounded curvature, implying
$\Vert J_h\Vert _{\infty}:= \sup_{(y,t) \in[0,1]^d} |J_h((y,t))| < \infty$.

Given $y \in\M$, let $N_y$ be the set of points in $[0,1]^d$ with
parameterization $(y,t)$ for some $t \in\R$. Define $T_y:= \{t \in\R
\dvtx (y,t) \in N_y \}$. This gives
\begin{eqnarray*}
&&\la^{-(d-1)/d} \E H^\xi(\P_\la,\M) \\
&&\qquad=
\la^{1/d} \int_{y \in\M} \int_{t \in T_y}
\E\xi_\la \bigl((y,t), \P_\la,\M \bigr) \bigl|J_h
\bigl((y,t) \bigr)\bigr| \ka \bigl((y,t) \bigr)\,dt\,dy.
\end{eqnarray*}
Let $t = \la^{-1/d}u$ to obtain
%
\begin{eqnarray}
\label{4a}&& \la^{-(d-1)/d} \E H^\xi(\P_\la,\M)\nonumber
\\
&&\qquad= \int_{y \in\M}
\int_{u \in\la^{1/d} T_y} \E\xi_\la \bigl( \bigl(y,
\la^{-1/d}u \bigr), \P_\la,\M \bigr) \bigl|J_h \bigl(
\bigl(y,\la^{-1/d}u \bigr) \bigr)\bigr|\\
&&\hspace*{98pt}{}\times \ka \bigl( \bigl(y,\la^{-1/d}u
\bigr) \bigr)\,du \,dy.\nonumber
\end{eqnarray}
By Lemma
\ref{L1}, for almost all $y \in\M$ and $u \in\R$, we have
%
\begin{equation}
\label{4b}\lim_{\la\to\infty} \E\xi_\la \bigl(
\bigl(y, \la^{-1/d}u \bigr), \P_\la,\M \bigr) = \E\xi \bigl((
\0_y,u), \H_{\ka(y)}, \HH_y \bigr).
\end{equation}
By
\eqref{mom}, for $y \in\M$, $u \in\R$, and $\la\in(0, \infty)$, the
integrand in \eqref{4a} is bounded by $G^{\xi,1}(|u|) \Vert J_h\Vert _\infty
\Vert \ka\Vert _\infty$, which is integrable with respect to the measure $du
\,dy$. 
Therefore, by the dominated convergence
theorem, the limit $\la^{1/d} T_y \uparrow\R$, the continuity of
$\ka$, and \eqref{4b}, we obtain \eqref{one}, namely
\[
\lim_{\la\to\infty} \la^{-(d-1)/d} \E H^\xi(
\P_\la,\M) = \int_{y
\in\M} \int_{-\infty}^{\infty}
\E \bigl[\xi \bigl((\0_y,u), \H_{\ka(y)}, \HH_y
\bigr) \bigr] \,du \,\ka(y) \,dy.
\]
To show \eqref{two}, we note
\begin{eqnarray*}
&&\la^{-2(d-1)/d} \E H^\xi(\P_\la,\M)^2\\
&&\qquad =
\la^{-2(d-1)/d}\\
&&\quad\qquad{}\times \biggl [ \la \int_{[0,1]^d} \E \bigl[
\xi_\la(x, \P_\la, \M)^2 \bigr] \ka(x) \,dx
\\
&&\hspace*{15pt}\qquad\quad{}+ \la^2 \int_{[0,1]^d} \int_{[0,1]^d}
\E\xi_\la(x, \P_\la, \M) \xi_\la(w,
\P_\la, \M) \ka(x) \ka(w) \,dx \,dw\biggr].
\end{eqnarray*}
The first integral goes to zero, since $\sup_{\la> 0}
\la^{1/d} \int_{[0,1]^d} \E\xi_\la(x, \P_\la, \M)^2 \ka(x) \,dx$ is
bounded. The second integral simplifies to
\[
\la^{2/d} \int_{[0,1]^d} \int_{[0,1]^d}
\E\xi_\la(x, \P_\la, \M ) \xi_\la(w,
\P_\la, \M) \ka(x) \ka(w) \,dx \,dw.
\]
As $\la\to\infty$, this tends to $\mu(\xi, \M)^2$ by
independence, proving the asserted $L^2$ convergence
of Theorem \ref{WLLN}.

To prove $L^1$ convergence we follow a truncation argument similar to
that for the proof of Proposition 3.2
in \cite{PY4}. Given $K > 0$, we put
\[
\xi^K(x, \X, \M):= \min \bigl(\xi(x,\X, \M), K \bigr).
\]
Then $\xi^K$ is homogenously stabilizing and uniformly bounded and,
therefore, by the first part of this proof
we get
%
\begin{equation}
\label{conv-b} \lim_{\la\to\infty} \la^{-(d-1)/d}
H^{\xi^K}( \P_\la, \M) = \mu \bigl(\xi^K, \M
\bigr)\qquad \mbox{in } L^2.
\end{equation}
Also, following the arguments around \eqref{4a}, we have
\begin{eqnarray*}
&&\bigl|\la^{-(d-1)/d} \bigl( \E H^\xi(\P_\la, \M) - \E
H^{\xi^K}(\P_\la, \M) \bigr)\bigr| \\
&&\qquad
\leq\int_{y \in\M} \int_{u \in\la^{1/d} T_y} \E[\cdots]
\bigl|J_h \bigl( \bigl(y,\la^{-1/d}u \bigr) \bigr)\bigr| \ka \bigl(
\bigl(y, \la^{-1/d}u \bigr) \bigr)\,du\,dy,
\end{eqnarray*}
where $\E[\cdots]:= \E[| \xi_\la((y,\la^{-1/d}u), \P_\la,\M) -
\xi^K_\la((y,\la^{-1/d}u), \P_\la,\M)|]$.
This expected difference tends to zero as $K \to\infty$, because the
moments condition~\eqref{mom}
with $p > 1$ implies that $| \xi_\la((y,\la^{-1/d}u), \P_\la,\M)
- \xi_\la^K((y,\la^{-1/d}u),\break \P_\la,\M)|$
is uniformly integrable. By monotone convergence, $\mu(\xi^K, \M)
\to\mu(\xi, \M)$ as $K \to\infty$.
Thus, letting $K \to\infty$ in \eqref{conv-b} we get the desired
$L^1$ convergence.
\end{pf*}

\begin{pf*}{Proof of Theorem \ref{Var}} We have
\begin{eqnarray*}
\la^{-(d-1)/d} \Var H^\xi(\P_\la, \M) &=&
\la^{1/d} \int_{[0,1]^d} \E\xi^2_\la(x,
\P_\la, \M)\ka(x) \H^d(dx)
\\
&&{}+ \la^{1 + 1/d} \int_{x \in[0,1]^d} \int_{w \in[0,1]^d}
\{\cdots \} \ka(x)\ka(w) \,dx \,dw,
\end{eqnarray*}
where
\[
\{\cdots \}:= \E\xi_\la(x, \P_\la\cup w, \M)
\xi_\la(w, \P_\la \cup x, \M) - \E\xi_\la(x,
\P_\la, \M) \E\xi_\la(w, \P_\la, \M).
\]
For a fixed $(y,t) \in\M\times\R$,
parameterize points $x \in[0,1]^d$ by $x_y:= (z_y, s_y)$, where $z_y
\in\HH_y$ and $s_y \in\R$.
Given $(y,t) \in[0,1]^d$ and $z_y \in\HH_y$, let
$S_{z_y}:=S_{z_y,t}:= \{s_y \in\R\dvtx (y,t) + (z_y, s_y) \in[0,1]^d \}$ and
let $Z_y:= [0,1]^d \cap\HH_y$. We have
%
\begin{eqnarray}
\label{var11}&& \la^{-(d-1)/d} \Var \bigl[H^\xi(
\P_\la, \M ) \bigr] \nonumber\\
&&\qquad= \la ^{1/d} \int_{[0,1]^d}
\E \xi_\la(x, \P_\la,\M)^2 \ka(x)\,dx
\nonumber
\\[-8pt]
\\[-8pt]
\nonumber
&&\qquad\quad{}+ \la^{1 + 1/d} \int_{y \in\M} \int_{T_y}
\int_{Z_y} \int_{S_{z_y}} \{\cdots \}
\bigl|J_h \bigl((y,t) \bigr)\bigr|\\
&&\hspace*{155pt}{}\times \ka \bigl((y,t) \bigr) \ka \bigl((y,t) +
(z_y,s_y) \bigr) \,ds_y
\,dz_y \,dt \,dy,\nonumber
\end{eqnarray}
where
\begin{eqnarray*}
\{\cdots\}&:=& \E\xi_\la \bigl((y,t), \P_\la\cup(y,t) +
(z_y,s_y), \M \bigr) \xi _\la \bigl((y,t) +
(z_y,s_y), \P_\la\cup(y,t), \M \bigr)
\\
&&{}- \E\xi_\la \bigl((y,t), \P_\la, \M \bigr) \E
\xi_\la \bigl((y,t) + (z_y,s_y), \P
_\la, \M \bigr).
\end{eqnarray*}
As in the proof of Theorem \ref{WLLN}, the first integral in \eqref
{var11} converges to
%
\begin{equation}
\label{first} \int_{\M} \int_{-\infty}^{\infty}
\E\xi^2 \bigl((\0_y,u), \H_{\ka
(y)},
\HH_y \bigr) \,du\, \ka(y) \,dy.
\end{equation}
In the second integral in \eqref{var11}, we let $t = \la^{-1/d}u, s_y
= \la^{-1/d} s, z_y = \la^{-1/d}z$ so that
$dz = \la^{(d-1)/d} \,dz_y$. These substitutions transform the
multiplicative factor
\[
\bigl|J_h \bigl((y,t) \bigr)\bigr| \ka \bigl((y,t) \bigr) \ka \bigl((y,t) +
(z_y,s_y) \bigr)
\]
into
%
\begin{equation}
\label{factor1} \bigl|J_h \bigl( \bigl(y,\la^{-1/d}u \bigr)
\bigr)\bigr| \ka \bigl( \bigl(y,\la^{-1/d}u \bigr) \bigr) \ka \bigl( \bigl(y,
\la^{-1/d}u \bigr) + \bigl(\la^{-1/d}z,\la^{-1/d}s \bigr)
\bigr),\hspace*{-10pt}
\end{equation}
they transform
the differential $\la^{1 + 1/d}\,ds_y \,dz_y \,dt \,dy$ into $\,ds \,dz \,du
\,dy$,
and, lastly, they transform
[recalling $x_y = (z_y,s_y)$] the covariance term $\{\cdots\}$ into
%
\begin{eqnarray}
\label{factor2} \{\cdots\}'&:= &\E\xi_\la \bigl( \bigl(y,\la
^{-1/d}u \bigr), \P_\la \cup \bigl(y,\la^{-1/d}u
\bigr) + \la^{-1/d}x_y, \M \bigr)\nonumber\\
&&{}\times \xi_\la \bigl(
\bigl(y, \la^{-1/d}u \bigr) + \la^{-1/d}x_y,
\P_\la \cup \bigl(y,\la ^{-1/d}u \bigr), \M \bigr) %
\\
&&{}- \E\xi_\la \bigl( \bigl(y,
\la^{-1/d}u \bigr), \P_\la, \M \bigr) \E\xi_\la
\bigl( \bigl(y, \la^{-1/d}u \bigr) + \la^{-1/d}x_y,
\P_\la, \M \bigr).\nonumber
\end{eqnarray}
The factor at \eqref{factor1} is bounded by $\Vert J_h\Vert _{\infty} \Vert \ka
\Vert _{\infty}^2$
and converges to $\ka(y)^2$, as $\la\to\infty$. 
By Lemma \ref{L2}, for almost all $y \in\M$, the covariance
term $\{\cdots\}'$ at \eqref{factor2} converges to 
\[
c^\xi \bigl((\0_y,u), (\0_y,u) + (z,s),
\H_{\ka(y)}, \HH_y \bigr).
\]
By Lemma
\ref{L3} as well as \eqref{mom}, the factor $\{\cdots\}'$ is dominated
by an integrable function of $(y,u,x_y) \in\M\times
\R\times\R^d$. By dominated convergence, together with the set limits
$\la^{1/d} Z_y \uparrow\R^{d-1}$, $\la^{1/d} S_{z_y} \uparrow\R$, and
$\la^{1/d} T_y \uparrow\R$
the second integral converges to
%
\begin{eqnarray}
\label{second} &&\int_{\M} \int_{\R^{d-1}}
\int_{-\infty}^{\infty}\int_{-\infty}^{\infty}
c^\xi \bigl((\0_y,u),(\0_y,u) + (z,s);
\H_{\ka(y)}, \HH_y \bigr)
\nonumber
\\[-8pt]
\\[-8pt]
\nonumber
&&\hspace*{86pt}{}\times \ka(y)^2 \,du \,ds \,dz
\,dy,
\end{eqnarray}
which is finite.
Combining \eqref{first} and \eqref{second}, we obtain Theorem \ref
{Var}.
\end{pf*}

\section{Proof of Theorem \texorpdfstring{\protect\ref{CLT}}{1.3}} \label{sec5}

Put $T_\la:= H^\xi(\P_\la,\M), \M\in\MM(d)$. We shall first
prove that Theorem
\ref{CLT} holds when $T_\la$ is replaced by a version $T'_\la$ on
input concentrated near $\M$. 
To show asymptotic normality of $T'_\la$, we follow the set-up of
\cite{PY5}, which makes use of dependency graphs, 
allowing applicability of Stein's method. We show that $T'_\la$ is
close to $T_\la$, thus yielding Theorem \ref{CLT}. This goes as
follows.

Put $\rho_\la:= {\beta}\log\la, s_\la:= \rho_\la\la^{-1/d} =
{\beta}\log
\la\cdot\la^{-1/d}$, ${\beta}\in(0, \infty)$ a constant to be
determined. Consider the collection of cubes $Q$ of the form
$\prod_{i=1}^d[j_i s_\la, (j_i + 1)s_\la)$, with all $j_i \in\Z$,
such that $\int_Q \ka(x) \,dx >
0$. 
Further, consider only cubes $Q$ such that $d(Q, \M) < 2s_\la$,
where for Borel subsets $A$ and $B$ of $\R^d$, we put $d(A, B):=
\inf\{|x - y|\dvtx  x \in A,  y \in B\}$. Relabeling if necessary,
write the union of the cubes as $\Q:= \bigcup_{i=1}^W Q_i$, where $W:=
W(\la) = \Theta((s_\la^{-1})^{d-1})$, because $\H^{d-1}(\M)<
\infty$.

We have $\operatorname{card}(Q_i \cap\P_\la):= N_i:= N(\nu_i)$, where $N_i$
is an independent Poisson random variable with parameter
\[
\nu_i:= \la\int_{Q_i} \ka(x) \,dx \leq\Vert \ka
\Vert _{\infty} \rho_\la^{d}. 
\]
We may thus write $\P_\la\cap\bigcup_{i=1}^W Q_i = \bigcup_{i=1}^W \{
X_{ij} \}_{j=1}^{N_i}$, where
for $1 \leq i \leq W$, we have $X_{ij}$ are i.i.d. on $Q_i$ with
density
\[
\ka_i(\cdot):= \frac{ \ka(\cdot)} { \int_{Q_i} \ka(x)\,dx} {\mathbf{1}}(Q_i).
\]
%

Define
\[
\tT_\la:= \sum_{x\in\P_\la\cap\Q}
\xi_\la(x, \P_\la, \M).
\]
Then by definition of $W$, $N_i$
and $X_{ij}$, we may write
\[
\tT_\la= \sum_{i=1}^{W} \sum
_{j = 1 }^{N_i} \xi_\la(X_{ij},
\P_\la, \M).
\]
%

As in \cite{PY5}, it is useful to consider a version $T'_{\la}$ of
$\tT_\la$ which has more independence between summands. This goes as
follows. For all $1 \leq i \leq W$ and all $j =1,2,\ldots ,$ recalling
the definition \eqref{expo}, let $R_{ij}:=R^\xi(X_{ij}, \P_\la, \M
)$ denote the radius of stabilization of $\xi$ at $X_{ij}$
if $1 \leq j \leq N_i$ and otherwise let $R_{ij}$ be zero. Put
$E_{ij}:= \{ R_{ij} \leq\rho_\la\}$, let
%
\begin{equation}
\label{Ela} E_\la:= \bigcap_{i=1}^{W}
\bigcap_{j=1}^{\infty} E_{ij}
\end{equation}
and define
\[
T'_{\la}:= \sum_{i=1}^{W}
\sum_{j = 1 }^{N_i} \xi_\la(X_{ij},
\P_\la, \M){\mathbf{1}}(E_{ij}).
\]

For all $1 \leq i \leq W$, define
\[
S_i:= S_{Q_i}:= \bigl(\Var T'_{\la}
\bigr)^{-1/2} \sum_{j =1}^{N_i}
\xi_\la ( X_{ij}, \P_\la, \M) {\mathbf{1}}
(E_{ij}).
\]
Note that $S_i$ and $S_j$ are independent if $d(Q_i, Q_j) > 2
\la^{-1/d} \rho_\la$. Put
\[
S_\la:= \bigl(\Var T'_{\la}
\bigr)^{-1/2} \bigl(T'_{\la} - \E
T'_{\la} \bigr) = \sum_{i=1}^{W}
(S_i - \E S_i).
\]

We aim to show that $T'_\la$ closely approximates $T_\la$, but
first we show that $\tT_\la$ closely approximates $T_\la$.

\begin{lemm} \label{L5} Given $\M\in\MM(d)$, let $G^{\xi,2}:=G^{\xi, 2,\M}$
satisfy \eqref{mom} and \eqref{strongmom}. Choose ${\beta}\in(0,\infty)$ so that
%
\begin{equation}
\label{bchoice} {\beta}\limsup_{|u| \to\infty} |u|^{-1} \log
G^{\xi,2}\bigl(|u|\bigr) < -8.
\end{equation}
Then
%
\begin{equation}
\label{L5a} \Vert \tT_\la- T_\la\Vert _2
= O \bigl(\la^{-3} \bigr)
\end{equation}
and
%
\begin{equation}
\label{L5b}|\Var\tT_\la- \Var T_\la| = O \bigl(
\la^{-2} \bigr).
\end{equation}
\end{lemm}

\begin{pf} Writing $\tT_\la= T_\la+ (\tT_\la- T_\la)$
gives
\[
\Var\tT_\la= \Var T_\la+ \Var[\tT_\la-
T_\la] + 2 \Cov (T_\la, \tT_\la-
T_\la).
\]
Now
\begin{eqnarray*}
&&\Var[\tT_\la- T_\la] \\
&&\qquad\leq\Vert \tT_\la-
T_\la\Vert _2^2 = \E \biggl(\sum
_{x \in\P_\la\setminus\Q} \xi_\la(x, \P_\la, \M)
\biggr)^2
\\
&&\qquad= \la^2 \int_{[0,1]^d \setminus\Q} \int_{[0,1]^d \setminus Q}
\E \bigl[\xi_\la(x,\P_\la,\M) \xi_\la(y,
\P_\la,\M) \bigr]\ka(x) \ka(y) \,dx \,dy.
\end{eqnarray*}
If $x \in[0,1]^d \setminus\Q$, then $d(x, \M) \geq\beta\log\la
\cdot\la^{-1/d}$. Thus, by \eqref{mom} and \eqref{strongmom}, for
large $\la$ we have $ \E\xi_\la(x,\P_\la,\M)^2 \leq G^{\xi,2}(\beta\log\la) \leq\exp(-8\log\la)
= \la^{-8}$.
Applying the Cauchy--Schwarz inequality to $\E\xi_\la(x,\P_\la,\M)
\xi_\la(y,\P_\la,\M)$ with $ x,y \in[0,1]^d \setminus Q$, we
obtain
%
\begin{equation}
\label{V1} \Vert \tT_\la- T_\la\Vert
_2^2 = O \bigl(\la^{-6} \bigr)
\end{equation}
which gives \eqref{L5a}. Also, since $\Vert T_\la\Vert _2 = O( \la)$ and
$\Vert \tT_\la- T_\la\Vert _2 = O(\la^{-3})$, another application of the
Cauchy--Schwarz inequality
gives
%
\begin{equation}
\label{V2} \Cov(T_\la, \tT_\la- T_\la)
\leq \Vert T_\la\Vert _2 \Vert \tT_\la-
T_\la\Vert _2 = O \bigl(\la^{-2} \bigr).
\end{equation}
%
Combining \eqref{V1} and \eqref{V2} gives \eqref{L5b}.
\end{pf}


\begin{lemm} \label{L6} Assume that $\xi$ satisfies the moment
conditions \eqref{mom}
and \eqref{strongmom} for some $p > q, q \in(2,3]$.
For ${\beta}$ large, we have
%
\begin{equation}
\label{L6a}\bigl \Vert T_\la- T'_\la\bigr\Vert
_2 = O \bigl(\la^{-3} \bigr)
\end{equation}
and
%
\begin{equation}
\label{L6b}\bigl|\Var T_\la- \Var T'_\la\bigr| =
O \bigl( \la^{-2} \bigr).
\end{equation}
\end{lemm}

\begin{pf} We have $\Vert T_\la- T'_\la\Vert _2 \leq\Vert T_\la-
\tT_\la\Vert _2 + \Vert \tT_\la- T'_\la\Vert _2 = O(\la^{-3}) + \Vert \tT_\la-
T'_\la\Vert _2$, by Lemma \ref{L5}. 
Note that $|\tT_\la- T'_\la| = 0$ on $E_\la$, with $E_\la$ defined
at \eqref{Ela}. Choosing ${\beta}$ large
enough, we have $P[E_\la^c] = O(\la^{-D})$ for any $D > 0$. By the
analog of Lemma 4.3 of \cite{PY5}, and using condition \eqref{mom},
we get for $q \in(2,3]$ that $\Vert \tT_\la- T'_\la\Vert _q = O(\la)$.
This,
together with the H\"older inequality, gives $\Vert (\tT_\la-
T'_\la){\mathbf{1}}(E_\la^c)\Vert _2 = O(\la^{-3})$, whence \eqref{L6a}.

To show \eqref{L6b}, we note that by \eqref{L5b} and the triangle
inequality, it is enough to show $|\Var\tT_\la- \Var T'_\la| = O(
\la^{-2})$. However, this follows by writing
\[
\Var\tT_\la= \Var T'_\la+ \Var \bigl[
\tT_\la- T'_\la \bigr] + 2 \Cov
\bigl(T'_\la, \tT_\la- T'_\la
\bigr),
\]
noting $\Var[\tT_\la- T'_\la] \leq\Vert \tT_\la-
T'_\la\Vert _2 = O(\la^{-3})$, and then using $\Vert T'_\la\Vert _2 = O(\la)$
and the Cauchy--Schwarz inequality to bound $\Cov(T'_\la, \tT_\la-
T'_\la)$ by $O( \la^{-2})$.
\end{pf}

Now we are ready to prove Theorem \ref{CLT}. Since \eqref{dCLT}
trivially holds for large enough $\la$ when $\Var T_\la< 1$, we may
without loss of generality assume $\Var T_\la\geq1$.

As in \cite{PY5}, we define a dependency graph $G_\la:= ({\cal
V}_\la, {\cal E}_\la)$ for $\{S_i \}_{i=1}^V$. The set ${\cal
V}_\la$ consists of the cubes $Q_1,\ldots,Q_V$ and edges $(Q_i, Q_j)$
belong to ${\cal E}_\la$ iff $d(Q_i, Q_j) < 2 \la^{-1/d} \rho_\la$.
Using Stein's method in the context of dependency graphs, we adapt
the proof in \cite{PY5} to show the asymptotic normality of $S_\la$,
$\la\to\infty$, and then use this to show the asymptotic normality
of $T_\la, \la\to\infty$. In \cite{PY5}, we essentially replace
the term $V = \Theta( \la/(\log\la)^d)$ by the smaller term $W =
\Theta(\la^{(d-1)/d}/ (\log\la)^{d-1} )$, and instead of (4.16) and
(4.17) of \cite{PY5}, we use \eqref{L6a} and \eqref{L6b}. Note that
for $p > q, q \in(2,3]$, we have $\Vert S_i\Vert _q = O((
\Var[T'_\la])^{-1/2} \rho_\la^{d(p + 1)/p} )$. We sketch the
argument as follows.

Let $c$ denote a generic constant whose value may change at each
occurrence. Following Section~4.3 of \cite{PY5} verbatim up to
(4.18) gives, via Lemma 4.1 of \cite{PY5}, with $p > q$, $q \in
(2,3]$ and $\theta:= c(\Var[T'_\la])^{-1/2} \rho_\la^{d(p + 1)/p}$:
%
\begin{eqnarray}
\label{diffbound} &&\sup_{t \in\R} \bigl\llvert P[S_\la
\leq t] - \Phi(t) \bigr\rrvert\nonumber \\
&&\qquad\leq c W \theta^q \leq c
\la^{(d - 1)/d} \rho_\la^{-(d-1)} \bigl( \Var
T'_{\la} \bigr)^{-q/2} \rho_\la^{d(p + 1)q/p}
\\
&&\qquad\leq c \la^{(d - 1)/d} \bigl( \Var[T_\la] \bigr)^{-q/2}
\rho_\la^{dq + 1},\nonumber
\end{eqnarray}
where we use
$\Var[T'_{\la}] \geq
\Var[T_\la)]/2$, which follows (for $\la$ large) from
\eqref{L6b}.

Follow verbatim the discussion between (4.18)--(4.20) of \cite{PY5},
with 
$V(\la)$ there replaced by $W$. Recall that $q \in(2,3]$ with $p >
q$. Making use of \eqref{L6a}, this gives the analog of (4.20) of
\cite{PY5}. In other words, this gives a constant $c$ depending on
$d, \xi, p$, and $q$ such that for all $\la\geq2$ the inequality
\eqref{diffbound} becomes
%
\begin{eqnarray}
\label{**} &&\sup_{t \in\R} \bigl\llvert P \bigl[ \bigl(\Var
T'_\la \bigr)^{-1/2} (T_\la- \E
T_\la) \leq t \bigr] - \Phi(t) \bigr\rrvert
\nonumber
\\[-8pt]
\\[-8pt]
\nonumber
&&\qquad\leq c \la^{(d - 1)/d} ( \Var T_\la)^{-q/2}
\rho_\la^{dq + 1} + c\la^{-2}.
\end{eqnarray}
By \cite{BY2,PeEJP}, we have $\Var T_\la= O(\la)$ and so $c \la
^{-2}$ is
negligible with respect to the first term on the right-hand side of
\eqref{**}.







Finally we replace $\Var T'_\la$ by $\Var T_\la$ on the left-hand
side of \eqref{**}.
As in \cite{PY5}, we have by the triangle inequality
%
\begin{eqnarray}
\label{0513} &&\sup_{t \in\R} \bigl\llvert P \bigl[(\Var
T_{\la})^{-1/2} (T_{\la} - \E T_{\la}) \leq
t \bigr] - \Phi(t) \bigr\rrvert
\nonumber\\
&&\qquad\leq\sup_{t \in\R} \biggl\llvert P \biggl[ \bigl(\Var
T'_{\la} \bigr)^{-1/2} (T_{\la
} - \E
T_{\la}) \leq t \cdot \biggl( {\Var T_{\la} \over\Var T'_{\la} }
\biggr)^{1/2} \biggr]
\nonumber
\\[-8pt]
\\[-8pt]
\nonumber
&&\hspace*{139pt}\qquad{}- \Phi \biggl(t \biggl( {\Var T_{\la} \over\Var T'_{\la} }
\biggr)^{1/2} \biggr) \biggr\rrvert
\\
&&\qquad\quad{}+ \sup_{t \in\R} \biggl\llvert \Phi \biggl(t \biggl(
{\Var T_{\la} \over\Var
T'_{\la} } \biggr)^{1/2} \biggr) - \Phi(t) \biggr\rrvert.\nonumber
\end{eqnarray}

We have
\[
\biggl\llvert \sqrt{ {\Var T_\la\over\Var T'_\la} } - 1 \biggr\rrvert \leq
\biggl\llvert {\Var T_\la\over\Var T'_\la} - 1 \biggr\rrvert = O \bigl(
\la^{-2} \bigr).
\]
Let $\phi:= \Phi'$ be the density of $\Phi$. Following the analysis
after (4.21) of
\cite{PY5}, we get
\begin{eqnarray*}
\sup_{t \in\R} \biggl\llvert \Phi \biggl(t \sqrt{
{\Var T_\la\over\Var
T'_\la} } \biggr) - \Phi(t) \biggr\rrvert &\leq& c \sup
_{t \in\R} \biggl( \biggl({ |t | \over
\la^2 } \biggr) \Bigl(
\sup_{u \in[ t - tc/\la^2,  t + tc/\la^2]} \phi(u) \Bigr) \biggr)
\\
&=& O \bigl(\la^{-2}
\bigr).
\end{eqnarray*}
This gives \eqref{dCLT} as desired.

\section{Proofs of Theorems \texorpdfstring{\protect\ref{main1}--\protect\ref{main4}}{2.1--2.5}}
We first give a general result useful in
proving versions of Theorems
\ref{WLLN}--\ref{CLT} for binomial input. Say that $\xi$ is \emph{binomially
exponentially stabilizing} with respect to the pair $(\X_n, \M)$ if
for all $x \in\R^d$ there is a radius of stabilization $R:=R^\xi(x,
\X_n, \M) \in(0, \infty)$
a.s. such that
%
\begin{equation}
\label{bi-expo} \xi_n \bigl(x, \X_n \cap
B_{n^{-1/d}R}(x), \M \bigr) = \xi_n \bigl(x, \bigl(
\X_n \cap B_{n^{-1/d}R}(x) \bigr) \cup\A, \M \bigr)
\end{equation}
for all locally finite $\A\subset\R^d \setminus B_{n^{-1/d}R}(x)$,
and moreover, the tail\break probability
$\tilde{\tau}(t):= \tilde{\tau}(t, \M):= \sup_{n \geq1, x \in\R
^d} P[R(x, \X_n, \M) > t]$ satisfies\break $\limsup_{t \to\infty}
t^{-1} \log\tilde{\tau}(t) < 0$.

\begin{lemm} \label{dePoiss} Let $\M\in\MM(d)$.
Let $\xi$ be exponentially stabilizing \eqref{expo},
binomially exponentially stabilizing \eqref{bi-expo}, and
assume the moment conditions \eqref{mom} and~\eqref{strongmom} hold for some $p > 2$. 
If there is constant $c_1 \in(0, \infty)$ such that
%
\begin{equation}
\label{Bd1} P \bigl[\bigl |\xi_n(X_1, \X_n,
\M)\bigr| \geq c_1 \log n \bigr] = O \bigl(n^{- 1 - 2/(1-
1/p)} \bigr),
\end{equation}
and if $N(n)$ is an independent Poisson random variable with parameter
$n$, then
%
\begin{equation}
\label{Bd0} \bigl|\Var H^\xi(\X_{n}, \M) - \Var
H^\xi(\X_{N(n)}, \M) \bigr| = o \bigl(n^{(d-1)/d} \bigr).
\end{equation}
\end{lemm}

\begin{pf} Let $D:= 2/(1- 1/p)$. By \eqref{Bd1}, there
is an event $F_{n,1}$, with $P[F_{n,1}^c] = O(n^{-D})$ such that on
$F_{n,1}$ we have
%
\begin{equation}
\label{Bnd1} \max_{1 \leq i \leq n + 1} \bigl|\xi_n(X_i,
\X_n, \M)\bigr| \leq c_1 \log n.
\end{equation}
As in the proof of
Theorem \ref{CLT}, put $s_n:= \beta\log n/n^{1/d}$, $\Q:=\Q(n):=
\bigcup_{i=1}^W Q_i$, where $d(Q_i, \M) < 2s_n$, $W:= W(n) =
O((s_n^{-1})^{d-1})$, and ${\beta}$ is a constant to be determined.
Consider the event $F_{n,2}$ such that for all $1 \leq i \leq n + 1$,
we have $\xi_n(X_i, \X_n, \M) = \xi_n(X_i, \X_n \cap B_{s(n)}(X_i),
\M)$. By binomial exponential stabilization \eqref{bi-expo} and for
${\beta}$ large enough, we have $P[F_{n,2}^c] = O(n^{-D})$. Define for
all $n = 1,2,\ldots $
\[
\tT_n:= \sum_{X_i \in\X_n \cap\Q_n}
\xi_n(X_i, \X_n \cap\Q_n, \M).
\]
%
As in Lemma \ref{L5}, for ${\beta}$ large we have the generous bounds
\[
\bigl|\Var H^\xi(\X_n, \M) - \Var\tT_n \bigr| = o
\bigl(n^{(d-1)/d} \bigr)
\]
and
\[
\bigl|\Var H^\xi(\X_{N(n)}, \M) - \Var\tT_{N(n)} \bigr| = o
\bigl(n^{(d-1)/d} \bigr).
\]
Therefore, to show \eqref{Bd0}, it is enough to show
%
\begin{equation}
\label{Bd3} \bigl|\Var\tT_n - \Var\tT_{N(n)}\bigr| = o
\bigl(n^{(d-1)/d} \bigr).
\end{equation}
Write
$\xi_n(X_i, \X_n)$ for $\xi_n(X_i, \X_n, \M)$. If $X_i \in
B_{s_n}^c(X_{n + 1}), 1 \leq i \leq n$, then on $F_{n,2}$ we have
$\xi_n(X_i, \X_n) = \xi_n(X_i, \X_{n + 1})$. On $F_{n,2}$, we thus
have
\[
|\tT_n - \tT_{n+1} | \leq\xi_n(X_{n + 1},
\X_{n + 1}) + \sum_{1
\leq i \leq n\dvtx X_i \in B_{s_n}(X_{n + 1})} \bigl|
\xi_n(X_i, \X_n) - \xi_n(X_i,
\X_{n + 1})\bigr |.
\]
Given a constant $c_2 \in(0, \infty)$, define
\[
F_{n,3}:= \bigl\{\operatorname{card} \bigl\{X_n \cap
B_{s_n}(X_{n + 1}) \bigr\} \leq c_2 \log n \bigr\}.
\]
Choose $c_2$ large such that $P[F_{n,3}^c] = O(n^{-D})$. On $F_{n,1}
\cap F_{n,2} \cap F_{n,3}$ we have by~\eqref{Bnd1} 
$|\tT_n - \tT_{n+1} | = O((\log n)^2)$. We deduce there is a $c_3$
such that on $F_{n,1} \cap F_{n,2} \cap F_{n,3}$ and all integers $l
\in\{1,\ldots,n\}$
%
\begin{equation}
\label{Bd6} |\tT_n - \tT_{n+l} | \leq
c_3l(\log n)^2.
\end{equation}

To show \eqref{Bd3}, we shall show
%
\begin{equation}
\label{Bd7a} |\Var\tT_n - \Var\tT_{N(n)}| = O \bigl((\log
n)^{4} n^{1 - 3/2d} \bigr).
\end{equation}
To show \eqref{Bd7a}, write
%
\begin{equation}
\label{triangle} \quad\Var\tT_n = \Var \tT_{N(n)} + (\Var
\tT_n - \Var\tT_{N(n)}) + 2 \operatorname{cov}( \tT_{N(n)},
\tT_n - \tT_{N(n)}).\hspace*{-10pt}
\end{equation}
The proof of Theorem
\ref{Var} shows $\Var\tT_{N(n)} = O(n^{(d-1)/2d})$, yielding
\begin{eqnarray*}
\operatorname{cov}( \tT_{N(n)}, \tT_n - \tT_{N(n)}) &\leq&
\sqrt{\Var\tT_{N(n)} } \cdot \bigl\Vert( \tT_n -
\tT_{N(n)})\bigr\Vert _2
\\
&=& O \bigl(n^{(d-1)/2d} \bigl\Vert (\tT_n - \tT_{N(n)} )
\bigr\Vert _2 \bigr).
\end{eqnarray*}
It is thus enough to show
%
\begin{equation}
\label{1*} \bigl\Vert( \tT_n - \tT_{N(n)})\bigr\Vert
_2^2 = O \bigl( (\log n)^8 n^{1 - 2/d}
\bigr),
\end{equation}
since the last
two terms in \eqref{triangle} are then $O((\log n)^{4} n^{1 -
3/2d})$. Relabel the $X_i, i \geq1$, so that $\X_n \cap\Q_n =
\{X_1,\ldots,X_{B(n,s_n)} \}, \X_{N(n)} \cap\Q_n = \{X_1,\ldots,\break X_{N(n
\cdot s_n)} \}$.

Put $E_n:= \{ B(n,s_n) \neq N(n \cdot s_n) \}$. There is a coupling
of $B(n, s_n)$ and $N(n \cdot s_n)$ such that $P[E_n] \leq s_n$. By
definition of $E_n$,
\begin{eqnarray*}
&&\bigl\Vert( \tT_n - \tT_{N(n)})\bigr\Vert _2^2
\\
&&\qquad= \int\biggl| \sum_{X_i \in\X_n \cap\Q_n} \xi_n(X_i,
\X_n \cap\Q_n) - \sum_{X_i \in\X_{N(n)} \cap\Q_n}
\xi_n(X_i, \X_{N(n)} \cap \Q_n)\biggr|^2\\
&&\hspace*{43pt}{}\times{\mathbf{1}}(E_n) \,dP.
\end{eqnarray*}

Now $|B(n, s_n) - N(n \cdot s_n) | \leq c_4 \log n \sqrt{n s_n} $ on
an event $F_{n,4}$ with $P[F_{n,4}^c] = O(n^{-D})$. Let $F_n:=
\bigcap_{i=1}^4 F_{n,i}$ and note that $P[F_n^c] = O(n^{-D})$. By
\eqref{Bd6}, we have
\begin{eqnarray}
\label{Bd7} \qquad&&\int\biggl| \sum_{X_i \in\X_n
\cap\Q_n}
\xi_n(X_i, \X_n \cap\Q_n) - \sum
_{X_i \in\X_{N(n)}
\cap
\Q_n} \xi_n(X_i,
\X_{N(n)} \cap\Q_n)\biggr|^2\nonumber\\
&&\hspace*{8pt}{}\times {\mathbf{1}}(E_n) {
\mathbf{1}}(F_n) \,dP
\\
&&\qquad
\leq \bigl(c_3c_4 \log n \sqrt{n s_n} (
\log n)^2 \bigr)^2.\hspace*{-20pt}\nonumber
\end{eqnarray}
For random variables $U$ and $Y$, we have
$\Vert UY\Vert _2^2 \leq\Vert U\Vert _{2p}^2 \Vert Y\Vert _{2q}^2,   p^{-1} + q^{-1} = 1$,
giving
\begin{eqnarray}\label{Bd8}
\bigl\Vert (\tT_n - \tT_{N(n)} ) {\mathbf{1}} \bigl(F_n^c
\bigr)\bigr \Vert ^2_2 &= &\Vert \tT_n -
\tT_{N(n)}\Vert _{2p}^2 \bigl\Vert {\mathbf{1}}
\bigl(F_n^c \bigr) \bigr\Vert _{2q}^2
\nonumber
\\[-8pt]
\\[-8pt]
\nonumber
 &=& O \bigl(n^2 \bigr) \bigl(P \bigl[F_n^c
\bigr] \bigr)^{1/q} = O(1).
\end{eqnarray}
%
Combining \eqref{Bd7}--\eqref{Bd8} yields \eqref{1*}
as desired:
\begin{eqnarray*}
\bigl\Vert( \tT_n - \tT_{N(n)})\bigr\Vert _2^2
&= &O \biggl( (\log n)^6 n s_n \int{\mathbf
1}(E_n) {\mathbf{1}}(F_n) \,dP \biggr) + O(1)
\\
&=& O \bigl((\log n)^6 ns_n P[E_n] \bigr) +
O(1)
\\
&=&O \bigl((\log n)^6 ns_n^2 \bigr) = O
\bigl(( \log n)^8 n^{1-2/d} \bigr).
\end{eqnarray*}
\upqed\end{pf}

\begin{pf*}{Proof of Theorem \ref{main1}} Recalling the
definition of $\nu^-$ at \eqref{defp}, we have
%
\begin{equation}
\label{nuident} \la \bigl(\Vol(A_\la) - \Vol(A) \bigr) = \sum
_{x \in\P_\la} \nu^-_{\la}(x, \P_\la,
\partial A) = H^{\nu^-}(\P_\la, \partial A),
\end{equation}
where the
last equality follows from \eqref{lsMn}. Therefore,
\[
\la^{(d + 1)/d} \Var \bigl[\Vol(A_\la) - \Vol(A) \bigr] =
\la^{-(d-1)/d}\Var \bigl[H^{\nu^-}(\P_\la, \partial A)
\bigr].
\]
Likewise,
\[
\la\Vol(A \bigtriangleup A_\la) = \sum_{x \in\P_\la}
\nu^+_{\la}(x, \P_\la, \partial A) = H^{\nu
^+}(
\P_\la, \partial A).
\]
It is therefore enough to show that $\nu^-$ and $\nu^+$ satisfy the
conditions of Theorem~\ref{CLT}. We show this for
$\nu^-$; similar arguments apply for $\nu^+$. Write $\nu$ for
$\nu^-$ in all that follows.

As seen in Lemma 5.1 of \cite{PeBer}, when $\ka$ is bounded away
from $0$ and infinity, the functional $\tilde{\nu}(x, \X):=
\Vol(C(x, \X))$ is homogeneously stabilizing and exponentially
stabilizing with respect to $\P_\la$.
Identical arguments show that $\nu$ is homogeneously stabilizing and
exponentially stabilizing with respect to $(\P_\la,\partial A)$.
The arguments in \cite{PeBer} may be adapted to show that $\nu$
satisfies the $p$-moment condition~\eqref{mom}, and we provide the
details.
For all $y \in\partial A, z \in\R^d, u \in\R$, we have
%
\begin{eqnarray}
\label{numom} &&\bigl|\nu_\la \bigl( \bigl(y, \la^{-1/d}u \bigr)
, \P_\la\cup z, \partial A\bigr) \bigr| %
\nonumber\\
&&\qquad\leq\omega_d \operatorname{diam} \bigl[C \bigl(
\bigl( \la^{1/d}y, u \bigr), \la^{1/d}(\P_\la\cup z)
\bigr) \bigr]^d \\
&&\qquad\quad{}\times{\mathbf{1}} \bigl(C \bigl( \bigl(
\la^{1/d}y, u \bigr), \la^{1/d}(\P_\la\cup z) \bigr)
\cap\partial A \neq \varnothing \bigr),\nonumber
\end{eqnarray}
where $\omega_d:= \pi^{d/2}[\Gamma(1 + d/2)]^{-1}$ is the volume of
the $d$-dimensional unit
ball. When $\ka$ is bounded away from zero, the factor $\operatorname{diam}
[C((\la^{1/d}y, u), \la^{1/d}(\P_\la\cup z))]^d$ has finite moments
of all orders, uniformly in $y$ and $z$ \cite{MY}. It may be seen
that $\E[{\mathbf{1}} (C((\la^{1/d}y, u), \la^{1/d}(\P_\la\cup z))
\cap
\partial A \neq\varnothing)]$
decays exponentially fast in $u$, uniformly in $y$ and $z$ (see, e.g.,
Lemma 2.2 of \cite{MY}), giving
condition \eqref{mom}. The
Cauchy--Schwarz inequality gives exponential decay \eqref{strongmom}
for $\nu$.

Thus, $\nu:= \nu^-$ satisfies all
conditions of Theorem \ref{CLT} and, therefore,
recalling~\eqref{nuident}, the first part of
Theorem \ref{main1} follows. The second part of Theorem \ref{main1}
follows from identical arguments involving $\nu:= \nu^+$.
\end{pf*}

\begin{pf*}{Proof of Theorem \ref{main2}} As seen above, $\nu$ is
homogeneously and exponentially stabilizing with respect to $(\P_\la,
\partial A)$. It remains only to establish that $\nu$ is
well-approximated by $\P_\la$ input on half-spaces \eqref{lin} and
we may then deduce the second part of Theorem \ref{main2} from Theorem
\ref{Var}. This goes as follows.

Fix $\partial A \in\MM_2(d), y \in\partial A$.
Translating $y$ to the origin, letting $\P_\la$ denote a Poisson
point process on $[0,1]^d - y$, letting $\partial A$ denote
$\partial A - y$, and using rotation invariance of $\nu$, it is
enough to show for all $w \in\R^d$ that
\[
\lim_{\la\to\infty}\E\bigl| \nu \bigl(w, \la^{1/d}
\P_\la, \la^{1/d}\partial A \bigr) - \nu \bigl(w,
\la^{1/d}\P_\la, \R^{d-1} \bigr) \bigr| = 0.
\]
Without loss of generality, we
assume, locally around the origin, that $\partial A \subset
\R_{-}^{d-1}$.

Let $\tilde{C}(w,\la^{1/d}\P_\la)$ be the union of $C(w,\la
^{1/d}\P_\la)$ and the Voronoi cells adjacent to
$C(w,\la^{1/d}\P_\la)$ in the Voronoi tessellation of $\P_\la$.
Consider the event
%
\begin{equation}
\label{defE}E(\la,w):= \bigl\{ \operatorname{diam} \bigl[\tilde{C} \bigl(w,
\la^{1/d}\P_\la \bigr) \bigr] \leq{\beta}\log\la \bigr\}.
\end{equation}
For ${\beta}$
large, we have
$P[E(\la,w)^c] = O(\la^{-2})$ (see, e.g., Lemma 2.2 of \cite{MY}).
Note that $\nu(w, \la^{1/d}\P_\la, \la^{1/d}\partial A)$ and $\nu(w,
\la^{1/d}\P_\la, \R^{d-1})$ have finite second moments, uniformly
in $w \in\R^d$
and $\la\in(0, \infty)$.
By the Cauchy--Schwarz inequality, for large $\beta
\in(0, \infty)$, we have for all $w \in\R^d$, 
\[
\lim_{\la\to\infty}\E\bigl| \bigl(\nu \bigl(w, \la^{1/d}
\P_\la, \la ^{1/d}\partial A \bigr) - \nu \bigl(w,
\la^{1/d}\P_\la, \R^{d-1} \bigr) \bigr) {\mathbf{1}}
\bigl(E(\la,w)^c \bigr) \bigr| = 0.
\]
It is therefore enough to show for all $w \in\R^d$ that
%
\begin{equation}\qquad
\label{diff} \lim_{\la\to\infty} \E\bigl| \bigl(\nu \bigl(w,
\la^{1/d}\P_\la, \la^{1/d}\partial A \bigr) - \nu
\bigl(w, \la^{1/d}\P_\la, \R^{d-1} \bigr) \bigr) {
\mathbf{1}} \bigl(E(\la,w) \bigr) \bigr| = 0.
\end{equation}
We first assume $w \in\R_{-}^{d-1}$;
the arguments with $w \in\R_{+}^{d-1}$ are nearly identical.
Moreover, we may assume $w \in
\la^{1/d}A$ for $\la$ large. 
Consider the (possibly degenerate) solid
%
\begin{equation}
\label{Delta} \Delta_\la(w):= \Delta_\la(w,{\beta}):=
\bigl(\R^{d-1}_{-} \setminus \la^{1/d}A \bigr) \cap
B_{2{\beta}\log\la}(w).
\end{equation}
Since $\partial A$ is
$C^2$, the solid $\Delta_\la(w)$ has maximal ``height'' $o((\Vert w\Vert  +
2{\beta}\log\la) \la^{-1/d})$ with respect to the hyperplane
$\R^{d-1}$. It follows that
\[
\Vol \bigl(\Delta_\la(w) \bigr) = O \bigl(\bigl(\Vert w\Vert + 2{\beta}
\log \la\bigr) \la^{-1/d} (2{\beta}\log \la)^{d-1} \bigr) = O \bigl((
\log \la)^d \la^{-1/d} \bigr).
\]
On the event $E(\la,w)$, the difference of the volumes $C(w,
\la^{1/d}\P_\la) \cap\la^{1/d}A^c$ and $C(w, \la^{1/d}\P_\la)
\cap
\R_+^{d-1}$
is at most 
$\Vol(\Delta_\la(w))$. Thus, 
\begin{eqnarray*}
&&\E\bigl| \bigl(\nu \bigl(w, \la^{1/d}\P_\la, \la^{1/d}
\partial A \bigr) - \nu \bigl(w, \la^{1/d}\P_\la,
\R^{d-1} \bigr) \bigr) {\mathbf{1}} \bigl(E(\la,w) \bigr)\bigr |
\\
&&\qquad\leq\Vol \bigl(\Delta_\la(w) \bigr) = O \bigl((\log
\la)^d \la^{-1/d} \bigr),
\end{eqnarray*}
which gives \eqref{diff}, and thus the variance asymptotics follow.

We next prove the first part of Theorem \ref{main2}, namely $\Var
\Vol(A_\la) =\break  \Omega(\la^{-(d-1)/d})$.
By assumption, there
is a $C^1$ subset $\Gamma$ of $\partial A$, with $\H^{d-1}(\Gamma) >
0$. Recalling $A \subset[0,1]^d$, subdivide $[0,1]^d$ into cubes
of edge length $l(\la):= (\lfloor\la^{1/d} \rfloor)^{-1}$. The
number $L(\la)$ of cubes having nonempty intersection with
$\Gamma$ satisfies $L(\la) = \Omega(\la^{(d-1)/d} )$, as otherwise
the cubes would partition $\Gamma$ into $o(\la^{(d-1)/d} )$ sets,
each of
$\H^{d-1}$ measure $O( (\la^{-1/d})^{d-1}) $, giving
$\H^{d-1}(\Gamma) = o(1)$, a contradiction.

Find a subcollection $Q_1,\ldots,Q_M$ of the $L(\la)$ cubes such that
$d(Q_i, Q_j) \geq2 \sqrt{d} l(\la)$ for all $i, j \leq M$, and $M =
\Omega(\la^{(d-1)/d})$. Rotating and translating $Q_i, 1 \leq i \leq
M$, by a distance\vspace*{1pt} at most $(\sqrt{d}/2) l(\la)$, if necessary, we
obtain a collection $\tQ_1,\ldots,\tQ_M$ of disjoint cubes (with faces
not necessarily parallel to a coordinate plane) such that:
\begin{itemize}
\item$d(\tQ_i, \tQ_j) \geq\sqrt{d} l(\la)$ for all $i, j \leq M$,

\item$\Gamma$ contains the center of each $\tQ_i$, here denoted
$x_i, 1 \leq i \leq M$.
\end{itemize}

By the $C^1$ property, $\Gamma$ is well-approximated locally around
each $x_i$ by a hyperplane $\HH_i$ tangent
to $\Gamma$ at $x_i$. Making a further rotation of $Q_i$, if
necessary, we may assume that $\HH_i$ partitions $\tQ_i$ into
congruent rectangular solids.

Write $\nu$ for $\nu^{-}$. We now exhibit a configuration of Poisson
points $\P_\la$ which has strictly positive probability, for which
$\la^{(d-1)/d} \Vol(A_\la)$ has variability bounded away from zero,
uniform in $\la$. Let $\overrightarrow{\0 n_i}, n_i \in\R^d$, be
the unit normal to $\Gamma$ at $x_i$.
Let $\varepsilon:= \varepsilon(\la):= l(\la)/8$ and subdivide each
$\tQ_i$ into $8^d$ subcubes of edge length~$\varepsilon$. Recall that
$B_r(x)$ denotes the Euclidean ball centered at $x \in\R^d$ with
radius $r$. Consider cubes $\tQ_i, 1 \leq i \leq M$, having these
properties:

\begin{longlist}[(a)]
\item[(a)] the subcubes of $\tQ_i$ having a face on $\partial\tQ_i$,
called the ``boundary subcubes,'' each contain at least one point from
$\P_\la$,

\item[(b)] $\P_\la\cap[ B_{\varepsilon/20 }(x_i - {\varepsilon\over10} n_i)
\cup B_{\varepsilon/20 }(x_i + {\varepsilon\over10} n_i)]
$ consists of a singleton, say $w_i$, and

\item[(c)] $\P_\la$ puts no other points in $\tQ_i$.
\end{longlist}

Relabeling if necessary, let $I:= \{1,\ldots,K\}$ be the indices of
cubes $\tQ_i$ having properties (a)--(c). It is easily checked that
the probability a given $\tQ_i$ satisfies property (a) is strictly
positive, uniform in $\la$. This is also true for properties
(b)--(c), showing that
%
\begin{equation}
\label{EK} \E K = \Omega \bigl(\la^{(d-1)/d} \bigr).
\end{equation}

Without loss of generality, we may assume that $A$ contains
$B_{\varepsilon/20 }(x_i - {\varepsilon\over10} n_i)$
but that $A \cap B_{\varepsilon/20 }(x_i + {\varepsilon\over10} n_i) =
\varnothing$.
Abusing notation, let $\Q:= \bigcup_{i=1}^K \tQ_{i}$ and put $\Q^c:= [0,1]^d \setminus\Q$. Let $\F_\la$ be the sigma algebra
determined by the random set $I$, the positions of
points of $\P_\la$ in all boundary subcubes,
and the positions of points $\P_\la$ in $\Q^c$. Given $\F_\la$,
properties (a) and (c)
imply that
$\Vol(C(w_i, \P_\la)) = \Omega(\varepsilon^d)$. Simple
geometry shows that when $w_i \in B_{\varepsilon/20 }(x_i - {\varepsilon
\over10} n_i)$ we have $\Vol(C(w_i, \P_\la) \cap A^c) =
\Omega(\varepsilon^d)$, that is the contribution to $A_\la$
by the cell $C(w_i, \P_\la)$
is $\Omega(\varepsilon^d)$. On the other hand, when $w_i \in
B_{\varepsilon/20 }(x_i + {\varepsilon\over10} n_i)$, then there is no
contribution to $A_\la$.
Moreover, in either case, the volume contribution to $A_\la$ arising
from points of $\P_\la$ in the boundary subcubes is modified by
$o(\varepsilon^d)$
regardless of the position of $w_i$.
Conditional on $\F_\la$, and using that $w_i$ is equally likely to
belong to either ball,
it follows
that $\Vol(A_\la\cap\tQ_i)$ has variability $\Omega(\varepsilon
^{2d}) = \Omega(\la^{-2})$, uniformly in $i \in I$,
that is,
%
\begin{equation}
\label{VV} \Var \bigl[\Vol( A_\la\cap\tQ_i) |
\F_\la \bigr] = \Omega \bigl(\la^{-2} \bigr),\qquad i \in I.
\end{equation}
By the conditional variance formula,
\begin{eqnarray*}
\Var \bigl[\Vol(A_\la) \bigr] &=& \Var \bigl[ \E \bigl[
\Vol(A_\la) | \F_\la \bigr] \bigr] + \E \bigl[ \Var \bigl[
\Vol(A_\la) | \F_\la \bigr] \bigr]
\\
&\geq&\E \bigl[ \Var \bigl[ \Vol(A_\la)| \F_\la \bigr]
\bigr]
\\
&=&\E \bigl[\Var \bigl[\Vol(A_\la\cap\Q) + \Vol \bigl(A_\la
\cap\Q^c \bigr) | \F_\la \bigr] \bigr].
\end{eqnarray*}


Given $\F_\la$, the Poisson--Voronoi tessellation of $\P_\la$ admits
variability only inside~$\Q$, that is $\Vol(A_\la\cap\Q^c)$ is
constant. Thus,
\begin{eqnarray*}
\Var \bigl[\Vol(A_\la) \bigr]& \geq&\E \bigl[\Var \bigl[
\Vol(A_\la\cap\Q) | \F_\la \bigr] \bigr]
\\
&=& \E \biggl[\Var \biggl[ \sum_{i \in I} \Vol(
A_\la\cap\tQ_i) | \F_\la \biggr] \biggr] =
\E\biggl[ \sum_{i \in I} \Var \bigl[\Vol( A_\la\cap
\tQ_i) | \F_\la \bigr]\biggr],
\end{eqnarray*}
since, given $\F_\la$, $\Vol( A_\la\cap\tQ_i), i \in I$, are
independent. By \eqref{EK} and \eqref{VV}, we have
\[
\Var \bigl[ \Vol(A_\la) \bigr] \geq c_5
\la^{-2} \E[K] = \Omega \bigl(\la^{-(d +
1)/d} \bigr),
\]
concluding the proof of Theorem \ref{main2} when $\nu$ is
set to $\nu^-$.

To show $\Var[ \Vol(A \bigtriangleup A_\la) ] = \Omega(\la^{-(d +
1)/d})$, consider cubes $\tQ_i, 1 \leq i \leq M$, having these
properties:

\begin{longlist}[(a$'$)]
\item[(a$'$)] the ``boundary subcubes,'' each contain at least one point
from $\P_\la$,

\item[(b$'$)] $\P_\la\cap B_{\varepsilon/20 }(x_i - {\varepsilon\over10}
n_i)$ consists of a singleton, say $w_i$, and

\item[(c$'$)] $\P_\la\cap[ B_{\varepsilon/20 }(x_i + {\varepsilon\over10}
n_i) \cup B_{\varepsilon/20 }(x_i + {\varepsilon} n_i)
]$ consists of a singleton, say $z_i$,

\item[(d$'$)] $\P_\la$ puts no other points in $\tQ_i$.
\end{longlist}

Let $I':= \{1,\ldots,K' \}$ be the indices of cubes $\tQ_i$ having
properties (a$'$)--(d$'$). Let $\F_\la$ be as above, with $I$ replaced
by $I'$. It suffices to notice that on $\F_\la$, we have
\begin{eqnarray*}
\Vol(A \bigtriangleup A_\la){\mathbf{1}} \bigl(z_i \in
B_{\varepsilon/20 }(x_i + {\varepsilon} n_i) \bigr) &\geq&2
\Vol(A \bigtriangleup A_\la){\mathbf{1}} \bigl(z_i \in
B_{\varepsilon/20 }(x_i + {\varepsilon/10} n_i) \bigr)\\
& =&
\Omega \bigl(\la^{-2} \bigr).
\end{eqnarray*}
From this, we may deduce the analog of \eqref{VV}, namely
\[
\Var \bigl[\Vol \bigl( (A \bigtriangleup A_\la) \cap
\tQ_i \bigr) | \F_\la \bigr] = \Omega \bigl(
\la^{-2} \bigr),\qquad i \in I,
\]
and follow the above arguments nearly
verbatim. This concludes the proof when $\nu$ is set to $\nu^+$.
\end{pf*}

\begin{pf*}{Proof of Theorem \ref{main4ab}} For any
$\partial A$, we have $|\nu_n^{\pm}(X_i, \X_n, \partial A)| \leq
\Vol
(C(X_i, \X_n)) \leq\omega_d( \operatorname{diam}[C(n^{1/d}X_i, n^{1/d}
\X_n)])^d$.
Let $D = 2/(1 - 1/p)$. Modifications of Lemma 2.2 of \cite{MY} show
that with
probability at least $1 - n^{-D-1}$ we have $
(\operatorname{diam}[C(n^{1/d}X_i, n^{1/d} \X_n)])^d = O(\log n)$, that is
to say $\nu^{\pm}$ satisfies \eqref{Bd1}. The discussion in Section~6.3 of
\cite{PeEJP} shows that the functionals $\nu^{+}$ and
$\nu^{-}$ are binomially exponentially stabilizing as at
\eqref{bi-expo}.
Theorem \ref{main4ab} follows from
Lemma~\ref{dePoiss}, Theorems \ref{main1}--\ref{main2}, and
Corollary \ref{main3}.
\end{pf*}

\begin{pf*}{Proof of Theorem \ref{SA}} 
It suffices to show that the functional $\alpha$
defining the statistics \eqref{alphastat}
satisfies the
conditions of Theorems \ref{WLLN} and \ref{Var} and then apply~\eqref{supersimplemean} and \eqref{supersimplevar} to the statistic
\eqref{alphastat} to obtain \eqref{areaWLLN} and \eqref{areaVar},
respectively. To do this, we shall follow the proof that the volume
functional $\nu$ defined at \eqref{defp} satisfies these conditions.
The proof that $\alpha$ is homogeneously stabilizing and satisfies
the moment condition \eqref{mom} follows nearly verbatim the proof
that $\nu$ satisfies these conditions, where we only need to
replace the factor $\omega_d \operatorname{diam} [C((\la^{1/d}y, u),
\la^{1/d}(\P_\la\cup z))]^d $ in \eqref{numom} by $\omega_{d-1}
\operatorname{diam} [C((\la^{1/d}y, u),  \la^{1/d}(\P_\la\cup z))]^{d-1}$.

To show that $\alpha$ is well-approximated by $\P_\la$ input on
half-spaces \eqref{lin}, by moment bounds on $\alpha$ and the
Cauchy--Schwarz inequality, it is enough to show the analog of
\eqref{diff}, namely for all $w \in\R^d$ that
%
\begin{equation}
\label{adiff} \hspace*{6pt}\qquad\lim_{\la\to\infty} \E\bigl| \bigl(\alpha \bigl(w,
\la^{1/d}\P_\la, \la^{1/d}\partial A \bigr) - \alpha
\bigl(w, \la^{1/d}\P_\la, \R^{d-1} \bigr) \bigr) {
\mathbf{1}} \bigl(E(\la,w) \bigr)\bigr | = 0,
\end{equation}
where $E(\la,w)$ is at \eqref{defE}.
Recalling the definition of $\Delta_\la(w)$ at \eqref{Delta}, define
\[
E_0(\la, w):= \bigl\{ \la^{1/d}
\P_\la\cap\Delta_\la(w) = \varnothing \bigr\}.
\]
Since the intensity measure of $\la^{1/d} \P_\la$ is upper bounded
by $\Vert \ka\Vert _{\infty}$, we have
%
\begin{eqnarray}
\label{EE} P \bigl[E_0(\la, w)^c \bigr] &=& 1 - P
\bigl[E_0(\la, w) \bigr] \leq1 - \exp \bigl(- \Vert \ka\Vert
_{\infty} \Vol \bigl(\Delta_\la(w) \bigr) \bigr)
\nonumber
\\[-8pt]
\\[-8pt]
\nonumber
&\leq& 1 - \exp \bigl(- c_6 (\log\la)^d \la^{-1/d}
\bigr) = O \bigl((\log\la)^d \la^{-1/d} \bigr).
\end{eqnarray}

On the event $E(\la, w) \cap E_0(\la, w)$, the scores $\al(w,
\la^{1/d}\P_\la, \la^{1/d}\partial A)$ and $\al(w, \la^{1/d}\P
_\la,
\R^{d-1})$ coincide. Indeed, on this event it follows that
$f$ is face of the boundary cell $C(w,
\la^{1/d} \P_\la)$ of $\la^{1/d}A_\la$ iff $f$ is a face
of a boundary cell of the Poisson--Voronoi tessellation of
$\R_{-}^{d-1}$. [If $f$ is a face of the boundary cell $C(w,
\la^{1/d} \P_\la), w \in\la^{1/d}A$, then $f$ is also a face of
$C(z, \la^{1/d} \P_\la)$ for some $z \in\la^{1/d}A^c$. If
$\Delta_\la(w) = \varnothing$, then $z$ must belong to $\R_{+}^{d-1}$,
showing that $f$ is a face of a boundary cell of the Poisson--Voronoi
tessellation of $\R_{-}^{d-1}$. The reverse implication is shown
similarly.]

On the other hand, since
\[
\bigl\llvert\! \bigl\llvert \bigl(\al \bigl(w, \la^{1/d}
\P_\la, \la^{1/d}\partial A \bigr) - \al \bigl(w,
\la^{1/d}\P_\la, \R^{d-1} \bigr) \bigr) {\mathbf{1}}
\bigl(E(\la, w) \bigr) \bigr\rrvert\! \bigr\rrvert _2 = O(1),
\]
and since by \eqref{EE} we have $P[E_{0}(\la, w)^c]
= O((\log\la)^d \la^{-1/d})$, it follows by the Cauchy--Schwarz
inequality that as $\la\to\infty$,
%
\begin{eqnarray}
\label{diff-alpha-2} &&\E\bigl| \bigl(\al \bigl(w, \la^{1/d}\P_\la,
\la^{1/d}\partial A \bigr) - \al \bigl(w, \la^{1/d}
\P_\la, \R^{d-1} \bigr) \bigr)
\nonumber
\\[-8pt]
\\[-8pt]
\nonumber
&&\hspace*{97pt}{}\times {\mathbf{1}} \bigl(E(\la, w)
\bigr) { \mathbf {1}} \bigl(E_{0}(\la, w)^c \bigr)\bigr| \to0.
\end{eqnarray}
Therefore, \eqref{adiff} holds and so $\alpha$ is
well-approximated by $\P_\la$ input on half-spaces and $\alpha$
satisfies all conditions of Theorems \ref{WLLN} and \ref{Var}. This
proves statements \eqref{areaWLLN}--\eqref{areaVar}. Note that
\eqref{areaWLLN1} follows from \eqref{measLLN}, proving Theorem
\ref{SA}. To show these limits hold when Poisson input is replaced
by binomial input $\X_n$ we shall show that $\alpha$ satisfies the
conditions of Lemma \ref{dePoiss}. Notice that $|\alpha_n(X_1, \X_n,
\partial A)| \leq\H^{d-1}(C(X_1, \X_n)) = O( \operatorname
{diam}[C(n^{1/d}X_1,\break n^{1/d} \X_n)]^{d-1}) = O( (\log
n)^{(d-1)/d})$ with probability at least $1 - n^{-D-1}$, that is
$\alpha$ satisfies condition \eqref{Bd1}, where $D = 2/(1 - 1/p)$.
The arguments in Section~6.3 of \cite{PeEJP} may be modified to show that $\alpha$ is
binomially exponentially stabilizing as at \eqref{bi-expo} and,
therefore, by Lemma \ref{dePoiss}, the limits
\eqref{areaWLLN}--\eqref{areaWLLN1} hold for input $\X_n$, as
asserted in remark (i) following Theorem~\ref{SA}.
\end{pf*}

\begin{pf*}{Proof of Theorem \ref{main4}} Orient $\partial A$
so that points $(y,t) \in A$, have positive $t$ coordinate. Notice that
$\zeta$ satisfies the decay condition \eqref{strongmom} for all $p
\in[1, \infty)$. Indeed, for all $z \in\R^d \cup\varnothing, y \in
\partial A, u \in(-\infty, \infty)$, and $\la\in(0, \infty)$, we have
\[
\bigl|\zeta_\la \bigl( \bigl(y, \la^{-1/d}u \bigr),
\P_\la\cup z, \partial A \bigr)\bigr| \leq {\mathbf{1}} \bigl( \bigl(K \oplus
\bigl(y, \la^{-1/d}u \bigr) \bigr) \cap A \cap\P_\la=
\varnothing \bigr).
\]
Now
\[
P \bigl[ \bigl(K \oplus \bigl(y, \la^{-1/d}u \bigr) \bigr) \cap A \cap
\P_\la= \varnothing \bigr] = \exp \bigl(-\la\Vol \bigl( \bigl(K \oplus
\bigl(y, \la^{-1/d}u \bigr) \bigr) \cap A \bigr) \bigr)
\]
decays exponentially fast in $|u| \in(0, \infty)$, uniformly in $y
\in\partial A$ and $\la\in(0, \infty)$ and therefore \eqref
{strongmom} holds for all
$p \in[1, \infty)$.

To see that $\zeta$ is homogeneously stabilizing as at \eqref{hom},
we argue as follows. Without loss of generality, let $\0$ belong to the
half-space $H$ with hyperplane $\HH$, as otherwise $\zeta(\0, {\cal
H}_\tau, \HH) = 0$. Now
$\zeta(\0, {\cal H}_\tau, \HH)$ is insensitive to point
configurations outside $K \cap H$ and so $R^\zeta(\H_\tau,\HH):=
\diam(K \cap H)$ is a radius of stabilization for~$\zeta$.

To show exponential stabilization of $\zeta$ as at \eqref{expo}, we
argue similarly. 
By definition of maximality, $\zeta_\la(x, \P_\la, \partial A)$ is
insensitive to point configurations outside
$(K \oplus x) \cap A$. In other words,
$\zeta(\la^{1/d}x, \la^{1/d}\P_\la, \la^{1/d}\partial A)$ is
unaffected by point configurations
outside
\[
K_\la(x):= \bigl(K \oplus\la^{1/d}x \bigr) \cap
\la^{1/d}A.
\]

Let $R(x):= R^\zeta(x, \P_\la, \partial A)$
be the distance between $\la^{1/d}x$ and the nearest point in $K_\la
(x) \cap\la^{1/d}\P_\la$, if
there is such a point; otherwise let $R(x, \P_\la, \partial A)$
be the maximal distance between $\la^{1/d}x$ and $K_\la(x) \cap
\partial(\la^{1/d}A)$, denoted here by $D(\la^{1/d}x)$. By the
smoothness assumptions on
the boundary, it follows that $K_\la(x) \cap B_t(x)$ has volume at
least $c_7 t^d$ for all
$ 0 \leq t \leq D(\la^{1/d}x)$. It follows that uniformly in $x \in
\partial A$ and $\la> 0$
%
\begin{equation}
\label{expz} P \bigl[R(x) > t \bigr] \leq\exp \bigl(-c_7
t^d \bigr),\qquad 0 \leq t \leq D \bigl(\la^{1/d}x \bigr).
\end{equation}
For $t \in
[D(\la^{1/d}x), \infty)$,
this inequality holds trivially. Moreover, we claim that $R(x)$ is
a radius of stabilization for $\zeta$ at $x$. Indeed, if $R(x) \in
(0, D(\la^{1/d}x))$, then $x$ is not maximal, and so
\[
\zeta \bigl(x, \la^{1/d}\P_\la\cap B_R(x),
\la^{1/d}\partial A \bigr)= 0.
\]
Point configurations outside
$B_R(x)$ do not modify the score $\zeta$. If $R(x) \in
[D(\la^{1/d}x), \infty)$ then
\[
\zeta \bigl(x, \la^{1/d}\P_\la\cap B_R(x),
\la^{1/d}\partial A \bigr)= 1
\]
and point configurations outside
$B_R(x)$ do not modify $\zeta$, since maximality
of $x$ is preserved. Thus, $R(x):= R^\zeta(x, \P_\la, \partial A)$ is
a radius of stabilization for $\zeta$ at $x$, it decays
exponentially fast by \eqref{expz}, and \eqref{expo} holds.

It remains to show that $\zeta$ is well-approximated by $\P_\la$
input on half-spaces
\eqref{lin}. As with the Poisson--Voronoi functional, it is enough to
show the
convergence \eqref{diff}, with $\nu$ replaced by $\zeta$ there.
However, since $\zeta$ is either $0$ or $1$,
we have that \eqref{diff} is bounded by the probability of the event
that $\la^{1/d} \P_\la$
puts points in the region $\Delta_\la(w)$ defined at \eqref{Delta}.
However, this probability tends to
zero as $\la\to\infty$, since the complement probability satisfies
\[
\lim_{\la\to\infty} P \bigl[ \la^{1/d} \P_\la\cap
\Delta_\la(w) = \varnothing \bigr] = \lim_{\la\to\infty}
\exp \bigl(- \Vol \bigl( \Delta_\la(w) \bigr) \bigr) = 1.
\]
This gives the required analog of \eqref{diff} for $\zeta$ and so
$\zeta_\la$ satisfies \eqref{lin}, which was to be shown. Thus,
Theorem \ref{main4} holds for Poisson input $\P_\la$, where we note
$\sigma^2(\zeta, \partial A) \in(0, \infty)$ by Theorem 4.3 of
\cite{BX}.
Straightforward modifications of the above arguments show that $\zeta$
is binomially exponentially stabilizing as at \eqref{bi-expo}. Now
$|\zeta| \leq1$, so
$\zeta$ trivially satisfies \eqref{Bd1}. Therefore, by Lemma \ref
{dePoiss}, Theorem \ref{main4} holds
for binomial input $\X_n$.

This completes the proof of
Theorem \ref{main4}, save for showing \eqref{zetaLLN}.
First notice that
%
\begin{equation}
\label{zetalim} \mu(\zeta, \partial A) = \int_{\partial A} \int
_{0}^{\infty} \E\zeta \bigl((\0_y,u),
\H_{1}, \HH _y \bigr) \ka(y)^{(d - 1)/d} \,du \,dy,
\end{equation}
which follows from \eqref{dWLLN} and $ \E\zeta((\0_y,u),
\H_{\tau}, \HH_y) = \E\zeta((\0_y,u \tau^{1/d}), \H_{1}, \HH_y)$.

The limit \eqref{zetalim}
further simplifies as follows. In $d = 2$, we have for $y = (v, F(v))
\in\partial A$ and all $u \in(0,\infty)$ that
\[
\E\zeta \bigl((\0_y,u), \H_{1}, \HH_y
\bigr) = \exp \biggl( - {u^2\over2} {(1
+ F'(v)^2)\over|F'(v)|} \biggr),
\]
where we use that a right triangle with legs on the coordinate axes,
hypotenuse distant $u$ from the origin and
having slope $m \in(0, -\infty)$ has area $u^2(1 + m^2)/2|m|$. Put
$b:= {(1 + F'(v)^2)/ 2 |F'(v)|}$
and $z = u^2 b$. Then
\begin{eqnarray*}
\mu(\zeta, \partial A)&=&\int_{\partial A} \int_{0}^{\infty}
\E \zeta \bigl((\0_y,u), \H_{1}, \HH_y
\bigr) \,du\, \ka(y)^{1/2} \,dy
\\
&=& {1 \over2} \int_{v \in[0,1]} \int
_0^{\infty} \exp(-z) (bz)^{-1/2} \sqrt{1 +
F'(v)^2} \ka \bigl(v, F(v) \bigr)^{1/2} \,dz
\,dv
\\
&=& {1 \over2} \Gamma \biggl({1 \over2} \biggr) \int
_{v \in[0,1]} b^{-1/2} \sqrt {1 + F'(v)^2}
\ka \bigl(v, F(v) \bigr)^{1/2} \,dv
\\
&=& {1 \over2} \Gamma \biggl({1 \over2} \biggr) \int
_{v \in[0,1]} 2^{1/2} \bigl|F'(v)\bigr|^{1/2}
\ka \bigl(v, F(v) \bigr)^{1/2} \,dv
\\
&=& \biggl({\pi\over2} \biggr)^{1/2} \int
_0^1 \bigl|F'(v)\bigr|^{1/2} \ka
\bigl(v, F(v) \bigr)^{1/2} \,dv.
\end{eqnarray*}
%

More generally, in $d > 2$, 
assume that $F$ is continuously differentiable with partials which are
negative and bounded
away from $0$ and $-\infty$. Let $y \in\partial A$ be given by $y =
(v, F(v)), v \in D$, and put
$F_i:= \partial F/\partial v_i$. Then for $u \in(0,\infty)$ we have
\[
\E\zeta \bigl((\0_y,u), \H_{1}, \HH_y
\bigr) = \exp \biggl( \frac{ - u^d (1 +
\sum_{i=1}^{d-1}
F_i'(v)^2)^{d/2}} {
d! \llvert  \prod_{i=1}^{d-1} F_i(v) \rrvert ^{-1} } \biggr).
\]
Let $z = u^d b$, where $b:= {1 \over d!} (1 + \sum_{i=1}^{d-1}
F_i'(v)^2)^{d/2} \llvert \prod_{i=1}^{d-1} F_i(v) \rrvert ^{-1}$. This
yields
\begin{eqnarray*}
\mu(\zeta, \partial A)&:=&\int_{\partial A} \int_{0}^{\infty}
\E \zeta \bigl((\0_y,u), \H_{1}, \HH_y
\bigr) \,du\, \ka(y)^{(d-1)/d} \,dy
\\
&= &({d!})^{1/d} d^{-1} \Gamma \bigl(d^{-1} \bigr)
\int_D \Biggl\llvert \prod_{i=1}^{d-1}
F_i(v) \Biggr\rrvert ^{1/d} \ka \bigl(v, F(v)
\bigr)^{(d-1)/d} \,dv,
\end{eqnarray*}
that is to say
\eqref{zetaLLN} holds.\vadjust{\goodbreak}
\end{pf*}

\section*{Acknowledgements} It is a pleasure to thank Y.
Baryshnikov for conversations related to Section~\ref{maxsection}
and the proof of Theorem \ref{main4}. 
Christoph Th\"ale kindly clarified the validity of \eqref{param} and
pointed me to \cite{Fr,HLW}, whereas Matthias Schulte shared \cite
{Sch} before its publication.




\printaddresses

\end{document}